\documentclass[11pt]{article}
\usepackage{amssymb,amsmath,color,amsfonts,float,amscd,amsthm,bm,mathrsfs} 

 \allowdisplaybreaks

\catcode`!=11
\let\!int\int \def\int{\displaystyle\!int}
\let\!lim\lim \def\lim{\displaystyle\!lim}
\let\!sum\sum \def\sum{\displaystyle\!sum}
\let\!sup\sup \def\sup{\displaystyle\!sup}
\let\!inf\inf \def\inf{\displaystyle\!inf}
\let\!cap\cap \def\cap{\displaystyle\!cap}
\let\!max\max \def\max{\displaystyle\!max}
\let\!min\min \def\min{\displaystyle\!min}

\catcode`!=12

\let\oldsection\section
\renewcommand\section{\setcounter{equation}{0}\oldsection}

\allowdisplaybreaks

\newtheorem{defn}{Definition}[section]
\newtheorem{thm}{Theorem}[section]
\newtheorem{prop}{Proposition}[section]
\newtheorem{lem}{Lemma}[section]
\newtheorem{cor}{Corollary}[section]
\newtheorem{rem}{Remark}[section]

\begin{document}
\title{Global Solvability of 2D MHD Boundary Layer Equations in Analytic Function Spaces}
\author{
{\bf Shengxin Li}\\[1mm]
\small lishengxin@sjtu.edu.cn, \\[1mm]
\small School of Mathematical Sciences,\\[1mm]
\small Shanghai Jiao Tong University,
Shanghai 200240, P.R.China\\[1mm]
{\bf Feng Xie}\\[1mm]
\small tzxief@sjtu.edu.cn, \\[1mm]
\small School of Mathematical Sciences, and LSC-MOE,\\[1mm]
\small Shanghai Jiao Tong University,
Shanghai 200240, P.R.China\\[1mm]
}
\date{}
\maketitle

\noindent{\bf Abstract:} In this paper we are concerned with the global well-posedness of solutions to magnetohydrodynamics (MHD) boundary layer equations in analytic function spaces. When the initial data is a small perturbation around a selected profile, and such a profile is governed by an one dimensional heat equation with a source term, we establish the global in time existence and uniqueness of analytic solutions to the two dimensional (2D) MHD boundary layer equations. It is noted that the far-field state of velocity is not required to be small initially, but decays to zero as time tends to infinity with suitable decay rates. The whole analysis is divided into two parts: When the initial far-field states are not small, we reformulate the original problem into a small perturbation problem by extracting a suitable background profile which is governed by a heat equation, then  we prove a long-time existence of solutions, and the lower bound of lifespan of solutions is given in term of the small parameter of initial perturbation; After establishing the long-time existence, combining the decay properties of far-field state and the properties of solutions to the heat equation, it is shown that both the solution obtained in the first part and the background profile indeed satisfy some smallness requirements when time goes to the lower bound of lifespan of solutions. Then we can show the global existence of solutions after this moment.


\vskip 2mm


\vskip 2mm

\section{Introduction}
In this paper, we consider global-in-time solvability of the following two-dimensional MHD boundary layer equations in the domain $\{(x, y)|x\in\mathbb{R}, y\in\mathbb{R}_+\}$.
\begin{equation}\label{11}
\begin{cases}
\partial_tu_1+u_1\partial_xu_1+u_2\partial_yu_1+\partial_xp=b_1\partial_xb_1+b_2\partial_yb_1+\partial_y^2u_1,\\
\partial_tb_1+\partial_y(u_2b_1-u_1b_2)=\kappa\partial_y^2b_1,\\
\partial_xu_1+\partial_yu_2=0,\quad \partial_xb_1+\partial_yb_2=0,
\end{cases}
\end{equation}
where the initial and boundary conditions are given as follows.
\begin{equation}\label{12}
u_1(t, x, y)|_{t=0}=u_{1, 0}(x, y),\qquad b_1(t, x, y)|_{t=0}=b_{1, 0}(x, y),
\end{equation}
and
\begin{equation}\label{13}
\begin{cases}
u_1|_{y=0}=u_2|_{y=0}=0, \\
\partial_yb_1|_{y=0}=b_2|_{y=0}=0.
\end{cases}
\end{equation}
The far-field states are denoted by $(\bar{u}, \bar{b})$.
\begin{equation}\label{14}
\lim\limits_{y\to+\infty}u_1=\bar{u}(t, x),\qquad \lim\limits_{y\to+\infty}b_1=\bar{b}(t, x),
\end{equation}
which satisfy Bernoulli's law
\begin{equation*}
\begin{cases}
\partial_t\bar{u}+\bar{u}\partial_x\bar{u}+\partial_xp=\bar{b}\partial_x\bar{b},\\
\partial_t\bar{b}+\bar{u}\partial_x\bar{b}=\bar{b}\partial_x\bar{u}.
\end{cases}
\end{equation*}

The solvability and properties of solutions to the initial boundary value problem (\ref{11})-(\ref{14}) play a key role in the study of high Reynolds number limit for the solutions to two dimensional (2D) magneto-hydrodynamics (MHD) equations in half-plane, where the no-slip boundary conditions are imposed on velocity, and the perfect-conducting boundary conditions are given for magnetic field on the boundary. Recently,  Liu, Yang and the second author \cite{LC1} established the local-in-time well-posedness of solutions to initial boundary value problem (\ref{11})-(\ref{14}) in weighted Sobolev spaces without the monotonicity condition on the initial velocity. Moreover, under the assumption that the viscosity and resistivity coefficients of 2D MHD equations are of the same order with respect to a small parameter, and the initial tangential magnetic field on the boundary does not degenerate, the same authors \cite{LC2} justified the validity of the Prandtl boundary layer ansatz in $L^\infty$ sense in finite regularity function spaces. For analytic initial data, Xie and Yang \cite{XF} obtained a lower bound of order $\varepsilon^{-2^-}$ for the lifespan of the solutions to (\ref{11})-(\ref{14}) when the initial datum is a small perturbation of a shear flow analytically in the order of $\varepsilon$. In the same framework, Liu and Zhang \cite{LN} succeeded in proving the global existence and the large time behavior of solutions to (\ref{11})-(\ref{14}), where the strength of far-field states are also required to be small initially, at least for the far-field state of velocity.

Before proceeding, it is necessary to recall the background and some related results of study of Prandtl boundary layer. It is noted that the MHD boundary layer equations are changed into the classical Prandtl boundary layer equations by taking $b_1=b_2=0$ in (\ref{11}). Such a system of  equations was first proposed by Prandtl in his seminal work \cite{PL} in 1904. Under a monotone assumption on the tangential velocity in the normal direction, Oleinik \cite{O1} obtained the local existence of classical solutions by using the Crocco transformation. The result together with other related works were written in the classical book \cite{O2}.  Recently, this important result was reproved by two groups
independently by using energy methods in \cite{A} and \cite{MN}. When imposing an additional favorable pressure condition, Xin and Zhang obtained a global-in-time weak solution in \cite{XZ}.
When the monotonicity condition is violated, the well-posedness of solutions still holds true in analytic function spaces and Gevrey regularity class. Sammaritino and Caflisch \cite{CR,S} established the local well-posedness theory of the Prandtl equations in three dimensional spaces and also justified the Prandtl boundary layer ansatz. Later, the analyticity requirement in the normal variable $y$ was removed by Lombardo, Cannone and Sammartino in \cite{LM} due to the partial viscous effect in the normal direction. Zhang and Zhang obtained the lifespan of small analytic solution to the classical Prandtl equations with small analytic data in \cite{Z}. Very recently, Paicu and Zhang \cite{PM1} proved the global existence and the large time decay estimate of solutions to Prandtl system with small initial data by using the weighted energy estimates. One also refer to \cite{I,K1,K2} for the well-posedness of Prandtl equations in analytic framework and \cite{GVMN,LW1,LW2} in Gevrey regularity framework.  It it noticed that most  well-posedness results of Prandtl equations and MHD boundary layer equations are concerning the local-in-time well-posedness except \cite{LN, PM1,XZ,Z} and other few works. Consequently, it is meaningful to study the global existence of smooth solution with a class of ``large" initial data, at least for the small perturbation around a ``large" profile. And this is the main motivation of study in this paper.

Below, we will introduce the main results and strategy of proof in this paper.  The global existence of solutions to (\ref{11})-(\ref{14}) in analytic function space
is the main goal of this paper, where the strength of far-field states are not required to be small initially comparing with the case considered in \cite{LN}.
The whole analysis in this paper can be divided into two parts: When the initial far-field states are not small, we reformulate the original problem into a small perturbation framework by extracting a suitable background profile, where this background solution is governed by an initial boundary value problem of a heat equation with a source term. Then we can achieve a long-time existence of solutions to the small perturbation problem, and give the lower bound of lifespan of solutions in term of the small parameter of initial perturbation. After establishing the long-time existence, combining the decay properties of the far-field state and the properties of solutions to the heat equation, it is shown that both the solution obtained in the first part and the background profile satisfy the smallness requirements in \cite{LN} (also \cite{PM1}) when time goes to the moment of lower bound of lifespan of solution. Then we can prove the global existence of solutions after this moment. Moreover, compared with the case in \cite{LN}, the decay rate requirements of the far-field state are also relaxed in this paper, see (\ref{112}) for details. However, it is necessary to point out that the far-field states considered in \cite{LN} also depend on $x$ variable.

Here and what follows, we only focus on the case that $\bar{b}=1$, $\kappa=1$ and $\bar{u}(t, x)=f(t)$ in the system $(\ref{11})$ and hence by Bernoulli's law the pressure term $\partial_xp=-f'(t)$. Then we choose a cut-off function $\varrho\in C^\infty [0, \infty)$ with
\begin{align}
\label{CF}
\varrho(y)=
\begin{cases}
1\quad \mbox{if} ~ y\ge 2,\\
0\quad \mbox{if}  ~ y\le 1.
\end{cases}
\end{align}
Denote $U_1=u_1-f(t)\varrho(y)$ and $(U_1, u_2, b_1, b_2)$ solves
\begin{equation}\label{15}
\begin{cases}
\partial_tU_1+(U_1+f(t)\varrho(y))\partial_xU_1+u_2\partial_y(U_1+f(t)\varrho(y))-b_1\partial_xb_1-b_2\partial_yb_1-\partial_y^2U_1=m(t, y),\\
\partial_tb_1+\partial_y(u_2b_1-(U_1+f(t)\varrho(y))b_2)=\partial_y^2b_1,\\
\partial_xU_1+\partial_yu_2=0,\qquad \partial_xb_1+\partial_yb_2=0,
\end{cases}
\end{equation}
where $m(t, y)=(1-\varrho(y))f'(t)+f(t)\varrho ''(y)$.

Note that a shear flow $(u^s(t, y), 0, 1, 0)$ is a trivial solution to the system $(\ref{15})$ with $u^s(t, y)$ solving
\begin{equation}\label{16}
\begin{cases}
\partial_tu^s(t, y)-\partial_y^2u^s(t, y)=m(t, y),\qquad (t, y)\in\mathbb{R}_+\times\mathbb{R}_+,\\
u^s(t, y)|_{y=0}=0,\quad\mbox{and}\quad\lim\limits_{y\to+\infty}u^s(t, y)=0,\\
u^s(t, y)|_{t=0}=0.
\end{cases}
\end{equation}

Then we rewrite the solution to (\ref{15}) as a perturbation $(u, v, b, g)$  of $(u^s(t, y), 0, 1, 0)$ by denoting
\begin{equation*}
\begin{cases}
U_1=u^s+u,\\
u_2=v,
\end{cases}
\qquad \qquad
\begin{cases}
b_1=1+b,\\
b_2=g.
\end{cases}
\end{equation*}
And (\ref{15}) becomes
\begin{equation}\label{17}
\begin{cases}
\partial_tu+(u+u^s+f(t)\varrho(y))\partial_xu+v\partial_y(u+u^s+f(t)\varrho(y))-(1+b)\partial_xb-g\partial_yb-\partial_y^2u=0,\\
\partial_tb+(u+u^s+f(t)\varrho(y))\partial_xb-g\partial_y(u+u^s+f(t)\varrho(y))-(1+b)\partial_xu+v\partial_yb-\partial_y^2b=0,\\
\partial_xu+\partial_yv=0,\qquad \partial_xb+\partial_yg=0,
\end{cases}
\end{equation}
with the initial data
\begin{equation}\label{18}
u(0, x, y)\triangleq u_0(x, y)=u_{1, 0}(x, y)-f(0)\varrho (y),\qquad b(0, x, y)\triangleq b_0(x, y)=b_{1, 0}(x, y)-1,
\end{equation}
the boundary conditions
\begin{equation}\label{19}
\begin{cases}
u|_{y=0}=0,\\
v|_{y=0}=0,
\end{cases}
\mbox{and}\qquad
\begin{cases}
\partial_yb|_{y=0}=0,\\
g|_{y=0}=0,
\end{cases}
\end{equation}
and the corresponding far-field conditions
\begin{equation}\label{110}
\lim\limits_{y\to+\infty}u=0,\qquad \lim\limits_{y\to+\infty}b=0.
\end{equation}

Since the divergence free conditions of $\partial_xu+\partial_yv=0$ and $\partial_xb+\partial_yg=0$, there exist two potential functions $(\varphi, \psi)$, such that $(u, b)=\partial_y(\varphi, \psi)$ and $(v, g)=-\partial_x(\varphi, \psi)$. Then by integrating the first two equations in (\ref{17}) with respect to $y$ variable over $[y, \infty)$ respectively, we obtain
\begin{equation}\label{111}
\begin{cases}
\partial_t\varphi+(u+u^s+f(t)\varrho(y)
)\partial_x\varphi+2\int_y^\infty\partial_y(u+u^s+f(t)\varrho(y))\cdot\partial_x\varphi\;dy\\
\qquad\qquad\qquad\qquad\qquad\qquad\qquad-(1+b)\partial_x\psi-2\int_y^\infty\partial_yb\cdot\partial_x\psi\;dy-\partial_y^2\varphi=0,\\
\partial_t\psi+(u+u^s+f(t)\varrho(y))\partial_x\psi-(1+b)\partial_x\varphi-\partial_y^2\psi=0,\\
\varphi |_{y=0}=\psi |_{y=0}=0\quad \mbox{and}\quad \lim\limits_{y\to+\infty}\varphi=\lim\limits_{y\to+\infty}\psi=0,\\
\varphi |_{t=0}=\varphi_0\triangleq-\int_y^\infty u_0\;dy^\prime,\quad \psi |_{t=0}=\psi_0\triangleq -\int_y^\infty b_0\;dy^\prime.
\end{cases}
\end{equation}
\begin{rem}
By divergence free conditions, one can show that $(u_0, b_0)$ satisfy the compatibility conditions of $\int_0^\infty u_0\;dy=\int_0^\infty b_0\;dy=0$ (see \cite{PM1}). Similarly, one has $\varphi|_{y=0}=\psi|_{y=0}=0$. These are exactly the boundary conditions of $\varphi$ and $\psi$ in (\ref{111}).
\end{rem}

The main results of this paper are stated as follows.
\begin{thm}
\label{T1.1}
Let $\delta>0$ and $f\in H^1(\mathbb{R}_+)$ which satisfies
\begin{equation}\label{112}
C_f\triangleq\int_0^\infty \langle t\rangle(|f(t)|+|f'(t)|)\;dt+(\int_0^\infty \langle t\rangle^{3}(f^2(t)+(f'(t))^2)\;dt)^{\frac12}<\infty
\end{equation}
with $\langle t\rangle=1+t$. And $u_0=\partial_y\varphi_0$, $b_0=\partial_y\psi_0$ satisfy the compatibility conditions $u_0(x, 0)=\partial_yb_0(x, 0)=0$ and $\int_0^\infty u_0\;dy=\int_0^\infty b_0\;dy=0$. Assume further that $G_0=u_0+\frac{y}{2\langle t\rangle}\varphi_0$ and $H_0=b_0+\frac{y}{2\langle t\rangle}\psi_0$ satisfy
\begin{equation}\label{113}
\|e^{\frac{y^2}8}e^{2\delta |D_x|}(u_0, b_0)\|_{\mathcal{B}^{\frac12, 0}}\le\varepsilon
\quad\mbox{and}\quad
\|e^{\frac{y^2}8}e^{\delta |D_x|}(G_0, H_0)\|_{\mathcal{B}^{\frac12, 0}}<\infty
\end{equation}
for some sufficiently small $\varepsilon$. Then the system (\ref{17})-(\ref{110}) admits a unique global-in-time solution $(u, b)$, which satisfies
\begin{equation}\label{114}
\|e^{\frac{y^2}{8\langle t\rangle}}e^{\frac\delta 2 |D_x|}(u, b)\|_{L^\infty(\mathbb{R}_+;\mathcal{B}^{\frac12, 0})}
+\|e^{\frac{y^2}{8\langle t\rangle}}e^{\frac\delta 2 |D_x|}\partial_y(u, b)\|_{L^2(\mathbb{R}_+;\mathcal{B}^{\frac12, 0})}
\le C\|e^{\frac{y^2}8}e^{2\delta |D_x|}(u_0, b_0)\|_{\mathcal{B}^{\frac12, 0}}.
\end{equation}
\end{thm}
The definition of anisotropic Besov spaces $\mathcal{B}^{\frac12, 0}$ will be given in section 2.
\begin{rem}
It should be emphasized that the smallness requirement in (\ref{113}) can be satisfied by a large class of initial data. Moreover, the spacial form of cutoff function $\varrho(y)$ in (\ref{CF}) is not essential.  Some other forms of cutoff functions also work for the whole analysis in this paper.
\end{rem}
\begin{rem}
It implies from (\ref{112}) that
\begin{equation}\label{115}
\begin{split}
|f(t)|&=|-\int_t^\infty f'(t')\;dt'|\le\int_t^\infty |\langle t'\rangle f'(t')|\frac1{\langle t'\rangle}\;dt\\
&\le\frac1{\langle t\rangle}\int_t^\infty \langle t'\rangle |f'(t')|\;dt\le\frac C{\langle t\rangle}
\end{split}
\end{equation}
for some positive constant $C$.
\end{rem}
\begin{rem}
The solution $u^s$ to (\ref{16}) has the following properties
\begin{equation}\label{116}
\|\partial_y^iu^s(t, \cdot)\|_{L_y^\infty}\le\frac C{\langle t\rangle^{i/2}},\quad i=1, 2,\qquad \|e^{\Psi(t, \cdot)}\partial_y^2u^s(t, \cdot)\|_{L_y^2}\le\frac C{\langle t\rangle^{3/4}},
\end{equation}
which can be proved by the properties of heat kernel, also refer to \cite{XF}. Here and after, $\Psi(t, y)\triangleq \frac {y^2}{8\langle t\rangle}$,
\end{rem}

Finally, the rest of this paper is organized as follows. In section 2, the function spaces used in this paper are defined and some technical lemmas are also given. Section 3  is devoted to establishing the lifespan of solution to the initial boundary value problem (\ref{17})-(\ref{110}). In section 4, the global existence of the solutions from the moment of lower bound of lifespan is proved based on the long-time existence obtained in section 3 and the decay properties of the far-field state.

Some notations used frequently in this paper are introduced below. We use the symbol $A \lesssim B$ to stand for $A \le CB$, where $C$ is a uniform constant which may vary from line to line. $(a, b)_{L_+^2}\triangleq\int_{\mathbb{R}_+^2}a(x, y)b(x, y)\;dxdy$ means the $L^2$ inner product of $a$, $b$ on $\mathbb{R}_+^2\triangleq\mathbb{R}\times\mathbb{R}_+$, and $(a, b)_{L_v^2}\triangleq\int_{\mathbb{R}_+}a(y)b(y)\;dy$. Finally, $L_T^p(L_h^q(L_v^r))$ denotes the space $L^p([0, T]; L^q(\mathbb{R}_x; L^r(\mathbb{R}_y^+)))$.

\section{Preliminaries}
In this paper, the Littlewood-Paley decomposition technique is used to define the function spaces, some basic facts about anisotropic Littlewood-Paley theory are collected in this section for the convenience of the readers. One also refer to \cite{BH} for more details.
\begin{equation}\label{21}
\Delta_k^ha=\mathcal{F}^{-1}(\phi(2^{-k}|\xi|)\hat{a}),\qquad S_k^ha=\mathcal{F}^{-1}(\chi(2^{-k}|\xi|)\hat{a}),
\end{equation}
where $\mathcal{F}$ and $\hat{a}$ always denote the partial Fourier transform of the distribution $a$ with respect to $x$ variable, that is, $\hat{a}(\xi, y)=\mathcal{F}_{x\to\xi}(a)(\xi, y)$. $\phi(\tau)$ and $\chi(\tau)$ are smooth functions such that
\begin{align*}
&Supp\, \phi \,\subset \{\tau\in\mathbb{R}\; |\; \frac34\le |\tau|\le\frac83\} \quad \mbox{and} \qquad\forall\tau>0,\ \sum\limits_{k\in\mathbb{Z}}\phi(2^{-k}\tau)=1,\\
&Supp\, \chi\, \subset \{\tau\in\mathbb{R} \;|\; |\tau|\le\frac83\} \qquad\quad \mbox{and} \qquad\chi(\tau)+\sum\limits_{k\ge0}\phi(2^{-k}\tau)=1.
\end{align*}
\begin{defn}
Let $s\in\mathbb{R}$. For $u\in\mathcal{S}_h^\prime(\mathbb{R}_+^2)$, which means that $u$ belongs to $\mathcal{S}^\prime(\mathbb{R}_+^2)$ and satisfies $\lim\limits_{k\to-\infty} \|\mathcal{S}_k^h u\|_{L^\infty}=0$, we set
\[
\|u\|_{\mathcal{B}^{s, 0}}\triangleq \|(2^{ks}\|\Delta_k^hu\|_{L_+^2})_{k\in\mathbb{Z}}\|_{\ell_1(\mathbb{Z})}.
\]
\end{defn}

In order to show a better description of the regularizing effect of the transport-diffusion equation, the following Chemin-Lerner type space $\tilde{L}_T^p(\mathcal{B}^{s, 0}(\mathbb{R}_+^2))$ is also needed.

\begin{defn}
Let $p\in [1, +\infty]$ and $T_0,\ T\in[0, +\infty]$. The space $\tilde{L}_T^p(T_0, T; \mathcal{B}^{s, 0}(\mathbb{R}_+^2))$ is defined as the completion of $C([T_0, T]; \mathcal{S}(\mathbb{R}_+^2)$ with the norm
\[
\|u\|_{\tilde{L}^p(T_0, T; \mathcal{B}^{s, 0})}\triangleq \sum\limits_{k\in\mathbb{Z}}2^{ks}\left(\int_{T_0}^T \|\Delta_k^hu\|_{L_+^2}^p\,dt\right)^{\frac1p}.
\]
If $p=\infty$, the above definition is modified by a usual change. When $T_0=0$, $\|u\|_{\tilde{L}_T^p(\mathcal{B}^{s, 0})}\triangleq \|u\|_{\tilde{L}^p(0, T; \mathcal{B}^{s, 0})}$ for simplicity of presentation.
\end{defn}

Furthermore, in order to use the Gronwall's type argument in the framework of Chemin-Lerner space $\tilde{L}_T^2(\mathcal{B}^{s, 0})$,
the time-weighted Chemin-Lerner type space was introduced in \cite{PM2}.

\begin{defn}
Suppose that $\mu(t)\in L_{loc}^1(\mathbb{R}_+)$ is a nonnegative function, and $t_0, t\in [0, +\infty]$. the time-weighted Chemin-Lerner type norm is defined as follows.
\begin{equation}\label{22}
\|u\|_{\tilde{L}_{t_0, t; \mu}^p(\mathcal{B}^{s, 0})}\triangleq\sum\limits_{k\in\mathbb{Z}}2^{ks}\left(\int_{t_0}^t \mu(t^\prime)\|\Delta_k^hu\|_{L_+^2}^p\,dt^\prime\right)^{\frac1p}.
\end{equation}
If $t_0=0$, we take the simplified notation of $\|u\|_{\tilde{L}_{0, t; \mu}^p(\mathcal{B}^{s, 0})}\triangleq \|u\|_{\tilde{L}_{t, \mu}^p(\mathcal{B}^{s, 0})}$.
\end{defn}

And the following anisotropic Bernstein estimates are also needed, one can refer to \cite{BJ} for details.
\begin{lem}
\label{L2}
Let $\mathcal{B}_h$ be a ball in $\mathbb{R}_h$, and $\mathcal{C}_h$ a ring in $\mathbb{R}_h$. Let $1\le p_2\le p_1\le\infty$ and $1\le q\le\infty$. Then it holds that\\
(1) If the support of $\hat{a}$ is included in $2^k\mathcal{B}_h$, then
\[
\|\partial_x^\ell a\|_{L_h^{p_1}(L_v^q)}\lesssim 2^{k(\ell+\frac1{p_2}-\frac1{p_1})}\|a\|_{L_h^{p_2}(L_v^q)}.
\]
(2) If the support of $\hat{a}$ is included in $2^k\mathcal{C}_h$, then
\[
\|a\|_{L_h^{p_1}(L_v^q)}\lesssim 2^{-k\ell}\|a\|_{L_h^{p_2}(L_v^q)}.
\]
\end{lem}

The Bony's decomposition in the horizontal variable $x$ is an effective tool to handle the estimates concerning the product of two distributions, see \cite{CJ1} for details.
\begin{equation}\label{23}
fg=T_f^hg+T_g^hf+R^h(f, g),
\end{equation}
where
\[
T_f^hg\triangleq\sum\limits_k \mathcal{S}_{k-1}^hf\Delta_k^hg, \quad R^h(f, g)\triangleq\sum\limits_k \tilde{\Delta}_k^hf\Delta_k^hg \quad \mbox{with} \quad \tilde{\Delta}_k^hf\triangleq\sum\limits_{|k-k^\prime |\le1} \triangle_{k^\prime}^hf.
\]

In addition, the following Poincar\'e type inequality, which is regarded as a special case of Treves inequality, can be found in \cite{I,LN,XF}.
\begin{lem}
\label{L2.1}
Let $\Psi(t, y)=\frac{y^2}{8\langle t\rangle}$ defined as above and $f$ is a smooth function defined on $\mathbb{R}_+^2$ such that $f|_{y=0}=0$ (or $\partial_yf|_{y=0}$) and $f|_{y=\infty}=0$. Then, it holds that
\begin{equation}\label{24}
\frac1{2\langle t\rangle}\|e^{\Psi}f\|_{L_v^2}^2\le\|e^{\Psi}\partial_yf
\|_{L_v^2}^2.
\end{equation}
\end{lem}
Below, we introduce some auxiliary functions which will be used later. First, the standard energy estimate argument to the system (\ref{16}) only yields that
\[
\|e^{\frac{y^2}{8\langle t\rangle}}\partial_yu^s(t)\|_{L_v^2}\le C\langle t\rangle^{-\frac34}.
\]
Yet, this estimate is not enough to guarantee that the quantity $\int_0^\infty \langle t\rangle^{\frac14}\|e^{\frac{y^2}{8\langle t\rangle}}\partial_yu^s(t)\|_{L_v^2}\;dt$ is finite, which is required for the later analysis.

To do this, we are going to introduce the primitive function $\varphi^s(t, y)$ of $u^s(t, y)$,  that is, $u^s(t, y)=\partial_y\varphi^s(t, y)$. It implies from (\ref{16}) that $\varphi^s$ satisfies the following initial boundary value problem.
\begin{equation}\label{25}
\begin{cases}
\partial_t\varphi^s(t, y)-\partial_y^2\varphi^s(t, y)=M(t, y),\qquad (t, y)\in\mathbb{R}_+\times\mathbb{R}_+,\\
\partial_y\varphi^s(t, y)|_{y=0}=0\quad\mbox{and}\quad\lim\limits_{y\to+\infty}\varphi^s(t, y)=0,\\
\varphi^s(t, y)|_{t=0}=0.
\end{cases}
\end{equation}
Here
\[
M(t, y)\triangleq -\int_y^\infty (1-\varrho(y'))\;dy'f'(t)+f(t)\varrho '(y).
\]
It is noted that $m(t, y)$ in (\ref{16}) equals to $\partial_yM$. Define
\begin{equation}\label{26}
G^s\triangleq u^s+\frac y{2\langle t\rangle}\varphi^s.
\end{equation}
By a direct calculation, the equation (\ref{25}) leads to
\[
\partial_t\left(\frac y{2\langle t\rangle}\varphi^s\right)-\partial_y^2\left(\frac y{2\langle t\rangle}\varphi^s\right)+\frac 1{\langle t\rangle}G^s=\frac y{2\langle t\rangle}M.
\]
Based on the definition (\ref{26}), combining with (\ref{16}), we arrive at
\begin{equation}\label{27}
\begin{cases}
\partial_tG^s-\partial_y^2G^s+\langle t\rangle^{-1}G^s=m+\frac y{2\langle t\rangle}M\triangleq \tilde{m}(t, y),\\
G^s|_{y=0}=0\quad\mbox{and}\quad\lim\limits_{y\to +\infty}G^s(t, y)=0,\\
G^s|_{t=0}=0.
\end{cases}
\end{equation}
In addition, in order to globally control the evolution of the analytic band of the solutions to (\ref{17}), the following two quantities are also needed.
\begin{equation}\label{28}
G\triangleq u+\frac y{2\langle t\rangle}\varphi,\qquad\qquad H\triangleq b+\frac y{2\langle t\rangle}\psi,
\end{equation}
which are inspired by the function $\mathfrak{g}\triangleq\partial_yu+\frac y{2\langle t\rangle}u$ introduced in \cite{I}. Where the weighted analytical norm of $\mathfrak{g}(t)$ decays like $\langle t\rangle^{-(\frac54)^-}$, this decay rate is faster than that of the weighted analytic norm of $u$ itself.
Then, it is helpful to consider the system of equations of $(G, H)$ as follows.
\begin{equation}\label{29}
\begin{cases}
&\partial_tG-\partial_y^2G+\frac G{\langle t\rangle}+(u+u^s+f(t)\varrho(y))\partial_xG-(1+b)\partial_xH+v\partial_yG-g\partial_yH\\
&\qquad+v\partial_y(u^s+f(t)\varrho(y))-\frac 1{2\langle t\rangle}v\partial_y(y\varphi)+\frac 1{2\langle t\rangle}g\partial_y(y\psi)\\
&\qquad-\frac y{\langle t\rangle}\int_y^\infty v\partial_y(u+u^s+f(t)\varrho(y))\;dy'
+\frac y{\langle t\rangle}\int_y^\infty g\partial_yb\;dy'=0,\\
&\partial_tH-\partial_y^2H+\frac H{\langle t\rangle}+(u+u^s+f(t)\varrho(y))\partial_xH-(1+b)\partial_xG+v\partial_yH-g\partial_yG\\
&\qquad-g\partial_y(u^s+f(t)\varrho(y))+\frac 1{2\langle t\rangle}g\partial_y(y\varphi)
-\frac 1{2\langle t\rangle}v\partial_y(y\psi)=0,\\
&G|_{y=0}=\partial_yH|_{y=0}=0,\quad \mbox{and}\quad \lim\limits_{y\to+\infty}G=\lim\limits_{y\to+\infty}H=0,\\
&G|_{t=0}=G_0,\qquad H|_{t=0}=H_0.
\end{cases}
\end{equation}
Next, we define the phase function $\Phi$ on $\mathbb{R}_+\times\mathbb{R}$ as in \cite{CJ1}.
\begin{equation}\label{210}
\Phi(t, \xi)\triangleq (\delta-\lambda\theta(t))|\xi|,
\end{equation}
and
\begin{equation}\label{211}
u_\Phi(t, x, y)=\mathcal{F}^{-1}_{\xi\to x}(e^{\Phi(t, \xi)}\hat{u}(t, \xi, y)),
\end{equation}
here, $\theta(t)$ is the key quantity to describe the evolution of the analytic band of $(u, b)$, the different definitions of $\theta(t)$ will be given in section 3 and section 4 respectively for different purpose.

We always assume that $t<T^*$ where $T^*$ is determined by
\begin{equation}\label{212}
T^*\triangleq\sup\{t>0, \theta(t)<\frac\delta\lambda\}.
\end{equation}
From (\ref{210}), for any $t<T^*$, it holds the following convex inequality
\begin{equation}\label{213}
\Phi(t, \xi)\le\Phi(t, \xi-\zeta)+\Phi(t, \zeta),\qquad \forall\xi, \zeta\in\mathbb{R}.
\end{equation}
The weighted function $\Psi(t, y)=\frac{y^2}{8\langle t\rangle}$ satisfies
\begin{equation}\label{214}
\partial_t\Psi(t, y)+2(\partial_y\Psi(t, y))^2=0.
\end{equation}

In the remainder of this section, we present the following lemma concerning the relationship between the function of $(G, H)$ defined in  (\ref{29}) and the solutions to (\ref{17}) and (\ref{111}), which will be frequently used in section 4.
\begin{lem}
\label{L3.2}
Let $(G, H)$ be defines by (\ref{29}) and let $(u, b)$ and $(\varphi, \psi)$ be smooth enough solutions to (\ref{17}) and (\ref{111}) respectively on $[0, T]$. Then for any $\gamma\in (0, 1)$ and $t\le T$, one has
\begin{align}
&\|e^{\gamma\Psi}\Delta_k^hu_\Phi(t)\|_{L_+^2}\lesssim\|e^\Psi\Delta_k^hG_\Phi(t)\|_{L_+^2},\quad \|e^{\gamma\Psi}\Delta_k^hb_\Phi(t)\|_{L_+^2}\lesssim\|e^\Psi\Delta_k^hH_\Phi(t)\|_{L_+^2};\label{215}\\
&\|e^{\gamma\Psi}\Delta_k^h\partial_yu_\Phi(t)\|_{L_+^2}\lesssim\|e^\Psi\Delta_k^h\partial_yG_\Phi(t)\|_{L_+^2},\quad \|e^{\gamma\Psi}\Delta_k^h\partial_yb_\Phi(t)\|_{L_+^2}\lesssim\|e^\Psi\Delta_k^h\partial_yH_\Phi(t)\|_{L_+^2};\label{216}\\
&\langle t\rangle^{-1}\|e^{\gamma\Psi}\Delta_k^h\partial_y(y\varphi)_\Phi(t)\|_{L_+^2}\lesssim\|e^\Psi\Delta_k^h\partial_yG_\Phi(t)\|_{L_+^2},\label{217a}\\
&\langle t\rangle^{-1}\|e^{\gamma\Psi}\Delta_k^h\partial_y(y\psi)_\Phi(t)\|_{L_+^2}\lesssim\|e^\Psi\Delta_k^h\partial_yH_\Phi(t)\|_{L_+^2}.\label{217b}
\end{align}
\end{lem}
The proof of this lemma can be found in \cite{LN, PM1} and we omit the detail here.

\section{Lifespan of Solutions}
Since the far-field state $f(t)$ are not required to be small initially in this paper, we can not apply the similar arguments as in \cite{LN,PM1} to  achieve the global existence of solutions directly. Consequently, we first establish the long-time existence of solutions to (\ref{17}), and give the lower bound of the lifespan of solutions. Then, by the decay properties of the far-field state $f(t)$, we will show that the far-field state $f(t)$ reach the requirement of smallness as in \cite{LN,PM1} when time goes to the lower bound of the lifespan of solutions. Moreover, the background profile defined in (\ref{16}) also become small at the moment of lower bound of the lifespan of solutions.

To this end, we introduce the following quantities.
\begin{equation}\label{41}
\tilde{u}=u-\partial_y(u^s+f(t)\varrho(y))\psi,\qquad \tilde{b}=b.
\end{equation}
Then $(\tilde{u}, \tilde{b})$ satisfies the following system of equations.
\begin{equation}\label{42}
\begin{cases}
\partial_t\tilde{u}-\partial_y^2\tilde{u}+(u+u^s+f(t)\varrho(y))\partial_x\tilde{u}+v\partial_y\tilde{u}-(1+b)\partial_x\tilde{b}-g\partial_y\tilde{b}\\
\qquad\qquad\qquad\qquad\qquad
-2\partial_y^2u^s\tilde{b}+v\partial_y^2u^s\psi-2f(t)\varrho ''(y)\tilde{b}+f(t)\varrho ''(y)v\psi=0,\\
\partial_t\tilde{b}-\partial_y^2\tilde{b}+(u+u^s+f(t)\varrho(y))\partial_x\tilde{b}+v\partial_y\tilde{b}-(1+b)\partial_x\tilde{u}-g\partial_y\tilde{u}\\
\qquad\qquad\qquad\qquad\qquad\qquad\qquad\qquad\qquad\qquad\,\,
-g\partial_y^2u^s\psi-f(t)\varrho ''(y)g\psi=0.
\end{cases}
\end{equation}
It is noted that the fact that $u^s$ is a solution to the heat equation is used. That is,
\[
\partial_tu^s-\partial_y^2u^s=m(t, y),\qquad \partial_t\partial_yu^s-\partial_y^3u^s=\partial_ym(t, y).
\]
By a direct calculation, the boundary conditions of $(\tilde{u}, \tilde{b})$ are given accordingly.
\begin{align}
&\tilde{u}|_{y=0}=0, \qquad\partial_y\tilde{b} |_{y=0}=0,\label{43}\\
&\lim\limits_{y\to\infty}\tilde{u}=0,\qquad\lim\limits_{y\to\infty}\tilde{b}=0. \label{44}
\end{align}
We then turn to show the existence of solution $(\tilde{u}, \tilde{b})$ to (\ref{42})-(\ref{44}) with the corresponding initial data,
\begin{equation}\label{45}
\tilde{u}(0, x, y)=u(0, x, y)-f(0)\varrho '(y)\psi_0,\qquad \tilde{b}(0, x, y)=b(0, x, y).
\end{equation}
Recall the definition of (\ref{210}), to avoid confusion, we use a new notation of $\eta(t)$ to describe the evolution of the analytic band to the solutions to (\ref{42})-(\ref{45}) and the corresponding phase function $\Phi(t, \xi)=(2\delta-\lambda\eta(t))|\xi|$ in this section, where
\begin{equation}\label{46}
\begin{cases}
\dot{\eta}(t)=\langle t\rangle^{\frac14}\left(\|e^\Psi\partial_y\tilde{u}_\Phi(t)\|_{\mathcal{B}^{\frac12, 0}}+\|e^\Psi\partial_y\tilde{b}_\Phi(t)\|_{\mathcal{B}^{\frac12, 0}}\right),\\
\eta |_{t=0}=0.
\end{cases}
\end{equation}
Note that
\begin{equation*}
\partial_yu=\partial_y\tilde{u}+\partial_y^2(u^s+f(t)\varrho(y))\psi+\partial_y(u^s+f(t)\varrho(y))\partial_y\psi
\quad\mbox{and}\quad
\tilde{b}=\partial_y\psi.
\end{equation*}
We can immediately get the following estimates by using (\ref{115}), (\ref{116}) and Lemma \ref{L2.1} that for any $s>0$,
\begin{equation}\label{47}
\|u\|_{\mathcal{B}^{s, 0}}\le\|\tilde{u}\|_{\mathcal{B}^{s, 0}}+C\|\tilde{b}\|_{\mathcal{B}^{s, 0}},\quad \|b\|_{\mathcal{B}^{s, 0}}=\|\tilde{b}\|_{\mathcal{B}^{s, 0}}.
\end{equation}
\begin{equation}\label{48}
\|\partial_yu\|_{\mathcal{B}^{s, 0}}\le\|\partial_y\tilde{u}\|_{\mathcal{B}^{s, 0}}+C\|\partial_y\tilde{b}\|_{\mathcal{B}^{s, 0}},\quad \|\partial_yb\|_{\mathcal{B}^{s, 0}}=\|\partial_y\tilde{b}\|_{\mathcal{B}^{s, 0}}.
\end{equation}
The main result in this section is stated as follows.
\begin{thm}
\label{T3.1}
Assume that $(u_0, b_0)$ satisfies (\ref{113}) and thus by (\ref{45}), it implies that $(\tilde{u}_0, \tilde{b}_0)$ satisfies
\[
\|e^{\frac{y^2}8}e^{2\delta |D_x|}(\tilde{u}_0, \tilde{b}_0)\|_{\mathcal{B}^{\frac12, 0}}\le C\varepsilon
\]
for some sufficiently small $\varepsilon$. Then there exists a positive time $T_\varepsilon$ with the size being greater than $\varepsilon^{-\frac32}$ so that (\ref{42})-(\ref{45}) admits a unique solution $(\tilde{u}, \tilde{b})$ which satisfies
\begin{equation*}
\int_0^t \langle t'\rangle^{\frac14}\|e^\Psi\partial_y(\tilde{u}, \tilde{b})_\Phi(t')\|_{\mathcal{B}^{\frac12, 0}}\;dt'\le C\langle t\rangle^{\frac23}\varepsilon
\end{equation*}
for any $t\in [0, T_\varepsilon]$, where $C\lambda\langle T_\varepsilon\rangle^{\frac23}\varepsilon=\delta$.
\end{thm}
\begin{rem}
Since the decay rate for the coefficient of linear term $-2\tilde{b}\partial_y^2(u^s+f(t)\varrho(y))$ in (\ref{42}) is not $L^1$ integrable with respect to time variable over $\mathbb{R}_+$, we can not achieve the global existence of solution to (\ref{42}) at this step. However, compared with the linear term $v\partial_y(u^s+f(t)\varrho(y))$ in (\ref{29}), the decay rate for the coefficient of  $-2\tilde{b}\partial_y^2(u^s+f(t)\varrho(y))$ is improved. Consequently, if the strength of $u^s$ is not required to be small, the lifespan obtained from (\ref{42}) should be larger than that from (\ref{29}).
\end{rem}
The proof of Theorem \ref{T3.1} relies on the following Proposition.
\begin{prop}\label{P3.1}
Let $(\tilde{u}, \tilde{b})$ be a smooth solution of (\ref{42})-(\ref{45}). Then for any nonnegative and non-decreasing function $\hbar(t)\in C^1(\mathbb{R}_+)$, there exists a suitably small positive constant $\alpha$ and a suitably large positive constant $K$, such that, for any $t<T^*$, it holds
\begin{equation}\label{49}
\begin{split}
&\|\hbar^{\frac12}e^\Psi(\tilde{u},\tilde{b})_\Phi\|_{\tilde{L}^\infty(0, t; \mathcal{B}^{\frac12, 0})}
+(1-2\sqrt{\alpha})\|\hbar^{\frac12}e^\Psi\partial_y\tilde{u}_\Phi\|_{\tilde{L}^2(0, t; \mathcal{B}^{\frac12,0})}
+(K-\frac {2C}{\sqrt{\alpha}})\|\hbar^{\frac12}e^\Psi\partial_y\tilde{b}_\Phi\|_{\tilde{L}^2(0, t; \mathcal{B}^{\frac12,0})}\\
&\le K\|\hbar^{\frac12}e^\Psi(\tilde{u},\tilde{b})_\Phi(0)\|_{\mathcal{B}^{\frac12, 0}}
+ \|\sqrt{\hbar '}e^\Psi\tilde{u}_\Phi\|_{\tilde{L}^2(0,t;\mathcal{B}^{\frac12,0})}
+ K\|\sqrt{\hbar '}e^\Psi\tilde{b}_\Phi\|_{\tilde{L}^2(0,t;\mathcal{B}^{\frac12,0})}.
\end{split}
\end{equation}
\end{prop}
Let us admit the above proposition for the time being and begin the proof of Theorem \ref{T3.1}.

Taking $\hbar(t)=1$ in (\ref{49}) and choosing suitable $\alpha$, $K>0$ such that both $1-2\sqrt{\alpha}$ and $K-\frac{2C}{\sqrt{\alpha}}$ are larger than zero give rise to
\begin{equation}\label{410}
\|e^\Psi(\tilde{u},\tilde{b})_\Phi\|_{\tilde{L}^\infty(0, t; \mathcal{B}^{\frac12, 0})}
+\|e^\Psi\partial_y(\tilde{u},\tilde{b})_\Phi\|_{\tilde{L}^2(0, t; \mathcal{B}^{\frac12,0})}\le
C\|e^{\frac{y^2}8}e^{2\delta |D_x|}(\tilde{u}_0,\tilde{b}_0)\|_{\mathcal{B}^{\frac12, 0}}.
\end{equation}
Next, it follows from Lemma \ref{L2.1} that
\[
\|\hbar^{\frac12}e^\Psi\partial_y(\tilde{u},\tilde{b})_\Phi\|_{\tilde{L}^2(0, t; \mathcal{B}^{\frac12,0})}\ge\frac1{\sqrt{2}}\|\langle t'\rangle^{-\frac12}\hbar^{\frac12}e^\Psi(\tilde{u},\tilde{b})_\Phi\|_{\tilde{L}^2(0, t; \mathcal{B}^{\frac12,0})}.
\]
And choosing $\hbar(t)=\langle t\rangle^{\frac16}$ in the above inequality and (\ref{49}), we obtain
\begin{align*}
\|\langle t\rangle^{\frac1{12}}e^\Psi&(\tilde{u},\tilde{b})_\Phi\|_{\tilde{L}^\infty(0, t; \mathcal{B}^{\frac12, 0})}
+\frac{1-2\sqrt{\alpha}}{\sqrt{2}}\|\langle t\rangle^{-\frac5{12}}e^\Psi\tilde{u}_\Phi\|_{\tilde{L}^2(0, t; \mathcal{B}^{\frac12,0})}
+\frac{K-\frac {2C}{\sqrt{\alpha}}}{\sqrt{2}}\|\langle t\rangle^{-\frac5{12}}e^\Psi\tilde{b}_\Phi\|_{\tilde{L}^2(0, t; \mathcal{B}^{\frac12,0})}\\
&\le K\|e^{\frac{y^2}8}(\tilde{u},\tilde{b})_\Phi(0)\|_{\mathcal{B}^{\frac12, 0}}
+ \frac1{\sqrt{6}}\|\langle t\rangle^{-\frac5{12}}e^\Psi\tilde{u}_\Phi\|_{\tilde{L}^2(0,t;\mathcal{B}^{\frac12,0})}
+ \frac{K}{\sqrt{6}}\|\langle t\rangle^{-\frac5{12}}e^\Psi\tilde{b}_\Phi\|_{\tilde{L}^2(0,t;\mathcal{B}^{\frac12,0})}.
\end{align*}
Let $\alpha$ small enough (for example $\alpha=\frac1{24}$), $K$ large enough, such that $\frac{1-2\sqrt{\alpha}}{\sqrt{2}}>\frac1{\sqrt{6}}$ and $\frac{K-\frac{2C}{\sqrt{\alpha}}}{\sqrt{2}}>\frac{K}{\sqrt{6}}$, we achieve that
\[
\|\langle t\rangle^{-\frac5{12}}e^\Psi(\tilde{u}, \tilde{b})_\Phi\|_{\tilde{L}^2(0,t;\mathcal{B}^{\frac12,0})}\le
C\|e^{\frac{y^2}8}e^{2\delta |D_x|}(\tilde{u}_0, \tilde{b}_0)\|_{\mathcal{B}^{\frac12, 0}},
\]
which immediately implies that
\begin{equation*}
\|\langle t\rangle^{\frac1{12}}e^\Psi(\tilde{u},\tilde{b})_\Phi\|_{\tilde{L}^\infty(0, t; \mathcal{B}^{\frac12, 0})}+\|\langle t\rangle^{\frac1{12}}e^\Psi\partial_y(\tilde{u},\tilde{b})_\Phi\|_{\tilde{L}^2(0, t; \mathcal{B}^{\frac12, 0})}\le C\|e^{\frac{y^2}8}e^{2\delta |D_x|}(\tilde{u}_0,\tilde{b}_0)\|_{\mathcal{B}^{\frac12, 0}}.
\end{equation*}
Therefore, in view of (\ref{46}),  it infers that
\begin{equation*}
\begin{split}
\eta(t)&\le\int_0^t \langle t'\rangle^{\frac14}\left(\|e^\Psi\partial_y\tilde{u}_\Phi(t')\|_{\mathcal{B}^{\frac12, 0}}+\|e^\Psi\partial_y\tilde{b}_\Phi(t')\|_{\mathcal{B}^{\frac12, 0}}\right)\;dt'\\
&\le\left(\int_0^t
(\|\langle t\rangle^{\frac1{12}}e^\Psi\partial_y\tilde{u}_\Phi(t')\|^2_{\mathcal{B}^{\frac12, 0}}
+\|\langle t\rangle^{\frac1{12}}e^\Psi\partial_y\tilde{b}_\Phi(t')\|^2_{\mathcal{B}^{\frac12, 0}})\;dt'\right)^{\frac12}\left(\int_0^t \langle t'\rangle^{\frac13}\;dt'\right)^{\frac12}\\
&\le C\langle t\rangle^{\frac23}\|e^{\frac{y^2}8}e^{2\delta |D_x|}(\tilde{u}_0, \tilde{b}_0)\|_{\mathcal{B}^{\frac12, 0}}.
\end{split}
\end{equation*}
In particular, under the assumption of Theorem \ref{T3.1}, we have
\begin{equation}\label{411}
\sup\limits_{t\in [0, T_\varepsilon^*]}\eta(t)\le\frac{\delta}{\lambda}\quad \mbox{for}\;T_\varepsilon^*\triangleq (\frac\delta {2\lambda C\varepsilon})^{\frac32}-1.
\end{equation}
This completes the proof of Theorem \ref{T3.1}.
\begin{rem}
The lower bound of lifespan in Theorem \ref{T3.1} can be improved to be $\varepsilon^{-2^-}$ as in \cite{XF}.
In fact, by taking $\hbar(t)=\langle t\rangle^{(\frac12)^-}$ in (\ref{49}), we can get the following estimate.
\[
\|\langle t\rangle^{(\frac14)^-}e^\Psi(\tilde{u},\tilde{b})_\Phi\|_{\tilde{L}^\infty(0, t; \mathcal{B}^{\frac12, 0})}+\|\langle t\rangle^{(\frac14)^-}e^\Psi\partial_y(\tilde{u},\tilde{b})_\Phi\|_{\tilde{L}^2(0, t; \mathcal{B}^{\frac12, 0})}\le C\|e^{\frac{y^2}8}e^{2\delta |D_x|}(\tilde{u}_0,\tilde{b}_0)\|_{\mathcal{B}^{\frac12, 0}},
\]
which ensures the lifespan is the order of $\varepsilon^{-2^-}$. However, the lower bound of $\varepsilon^{-\frac32}$ is enough for later use. In addition, it is noted that the method used here to achieve the lower bound of lifespan is different from that in \cite{XF}.
\end{rem}

Now, it suffices to prove Proposition \ref{P3.1}. In view of (\ref{42}) and (\ref{211}), we get
\begin{equation}\label{412}
\begin{cases}
\partial_t\tilde{u}_\Phi-\partial_y^2\tilde{u}_\Phi+\lambda\dot{\eta}(t)|D_h|\tilde{u}_\Phi+[(u+u^s+f(t)\varrho(y))\partial_x\tilde{u}]_\Phi+[v\partial_y\tilde{u}]_\Phi-
[(1+b)\partial_x\tilde{b}]_\Phi\\
 \qquad\qquad
 -[g\partial_y\tilde{b}]_\Phi-2\partial_y^2u^s\tilde{b}_\Phi+\partial_y^2u^s[v\psi]_\Phi-2f(t)\varrho ''(y)\tilde{b}_\Phi+f(t)\varrho ''(y)[v\psi]_\Phi=0,\\
\partial_t\tilde{b}_\Phi-\partial_y^2\tilde{b}_\Phi+\lambda\dot{\eta}(t)|D_h|\tilde{b}_\Phi+[(u+u^s+f(t)\varrho(y))\partial_x\tilde{b}]_\Phi+[v\partial_y\tilde{b}]_\Phi-
[(1+b)\partial_x\tilde{u}]_\Phi\\
\qquad\qquad\qquad\qquad\qquad\qquad\qquad\qquad
-[g\partial_y\tilde{u}]_\Phi+\partial_y^2u^s[g\psi]_\Phi-f(t)\varrho ''(y)[g\psi]_\Phi=0.
\end{cases}
\end{equation}
Here and what follows, we shall always denote $|D_h|$ to be the Fourier multiplier in the $x$ variable with symbol of $|\xi|$.

By applying the dyadic operator $\Delta_k^h$ to the above system and then taking the $L_+^2$ inner product of the resulting equations with $\hbar(t)e^{2\Psi}\Delta_k^h\tilde{u}_\Phi$ and $\hbar(t)e^{2\Psi}\Delta_k^h\tilde{b}_\Phi$ respectively, we get
\begin{equation}\label{413}
\begin{split}
\hbar(t)&(e^\Psi\Delta_k^h(\partial_t\tilde{u}_\Phi-\partial_y^2\tilde{u}_\Phi)|e^\Psi\Delta_k^h\tilde{u}_\Phi)_{L_+^2}
+\lambda\dot{\eta}(t)\hbar(t)(e^\Psi |D_h|\Delta_k^h\tilde{u}_\Phi |e^\Psi\Delta_k^h\tilde{u}_\Phi)_{L_+^2}\\
&+\hbar(t)(e^\Psi\Delta_k^h[(u+u^s+f(t)\varrho(y))\partial_x\tilde{u}]_\Phi |e^\Psi\Delta_k^h\tilde{u}_\Phi)_{L_+^2}
+\hbar(t)(e^\Psi\Delta_k^h[v\partial_y\tilde{u}]_\Phi |e^\Psi\Delta_k^h\tilde{u}_\Phi)_{L_+^2}\\
&-\hbar(t)(e^\Psi\Delta_k^h[(1+b)\partial_x\tilde{b}]_\Phi |e^\Psi\Delta_k^h\tilde{u}_\Phi)_{L_+^2}
-\hbar(t)(e^\Psi\Delta_k^h[g\partial_y\tilde{b}]_\Phi |e^\Psi\Delta_k^h\tilde{u}_\Phi)_{L_+^2}\\
&-\hbar(t)(2\partial_y^2u^se^\Psi\Delta_k^h\tilde{b}_\Phi |e^\Psi\Delta_k^h\tilde{u}_\Phi)_{L_+^2}
+\hbar(t)(\partial_y^2u^se^\Psi\Delta_k^h[v\psi]_\Phi |e^\Psi\Delta_k^h\tilde{u}_\Phi)_{L_+^2}\\
&-\hbar(t)(2f(t)\varrho ''(y)e^\Psi\Delta_k^h\tilde{b}_\Phi |e^\Psi\Delta_k^h\tilde{u}_\Phi)_{L_+^2}
+\hbar(t)(f(t)\varrho ''(y)e^\Psi\Delta_k^h[v\psi]_\Phi |e^\Psi\Delta_k^h\tilde{u}_\Phi)_{L_+^2}
=0,
\end{split}
\end{equation}
and
\begin{equation}\label{414}
\begin{split}
\hbar(t)&(e^\Psi\Delta_k^h(\partial_t\tilde{b}_\Phi-\partial_y^2\tilde{b}_\Phi)|e^\Psi\Delta_k^h\tilde{b}_\Phi)_{L_+^2}
+\lambda\dot{\eta}(t)\hbar(t)(e^\Psi |D_h|\Delta_k^h\tilde{b}_\Phi |e^\Psi\Delta_k^h\tilde{b}_\Phi)_{L_+^2}\\
&+\hbar(t)(e^\Psi\Delta_k^h[(u+u^s+f(t)\varrho(y))\partial_x\tilde{b}]_\Phi |e^\Psi\Delta_k^h\tilde{b}_\Phi)_{L_+^2}
+\hbar(t)(e^\Psi\Delta_k^h[v\partial_y\tilde{b}]_\Phi |e^\Psi\Delta_k^h\tilde{b}_\Phi)_{L_+^2}\\
&-\hbar(t)(e^\Psi\Delta_k^h[(1+b)\partial_x\tilde{u}]_\Phi |e^\Psi\Delta_k^h\tilde{b}_\Phi)_{L_+^2}
-\hbar(t)(e^\Psi\Delta_k^h[g\partial_y\tilde{u}]_\Phi |e^\Psi\Delta_k^h\tilde{b}_\Phi)_{L_+^2}\\
&+\hbar(t)(\partial_y^2u^se^\Psi\Delta_k^h[g\psi]_\Phi |e^\Psi\Delta_k^h\tilde{b}_\Phi)_{L_+^2}
-\hbar(t)(f(t)\varrho ''(y)e^\Psi\Delta_k^h[g\psi]_\Phi |e^\Psi\Delta_k^h\tilde{b}_\Phi)_{L_+^2}
=0.
\end{split}
\end{equation}
Let us integrate (\ref{413}) and ({\ref{414}) over  $[0,t]$, and handle term by term in the resulting equalities below.\\
$\bullet$ Estimate of $\int_0^t \hbar(t')(e^\Psi\Delta_k^h(\partial_t\tilde{u}_\Phi-\partial_y^2\tilde{u}_\Phi) |e^\Psi\Delta_h^h\tilde{u}_\Phi)_{L_+^2}\;dt^\prime$ and

\qquad\qquad$\int_0^t \hbar(t')(e^\Psi\Delta_k^h(\partial_t\tilde{b}_\Phi-\partial_y^2\tilde{b}_\Phi) |e^\Psi\Delta_h^h\tilde{b}_\Phi)_{L_+^2}\;dt^\prime$\\

First, by integration by parts, we have
\begin{equation}\label{415}
\begin{split}
\int_0^t& \hbar(t')(e^\Psi\Delta_k^h\partial_t\tilde{u}_\Phi |e^\Psi\Delta_k^h\tilde{u}_\Phi)_{L_+^2}\;dt^\prime\\
&=\frac12\|\hbar^{\frac12}e^\Psi\Delta_k^h\tilde{u}_\Phi(t)\|_{L_+^2}^2-\frac12\|\hbar^{\frac12}e^\Psi\Delta_k^h\tilde{u}_\Phi(0)\|_{L_+^2}^2\\
&-\frac12\int_0^t \hbar^\prime(t^\prime)\|e^\Psi\Delta_k^h\tilde{u}_\Phi(t^\prime)\|_{L_+^2}^2\;dt^\prime-\int_{0}^t\int_{\mathbb{R}_+^2} \hbar\partial_t\Psi |e^\Psi\Delta_k^h\tilde{u}_\Phi|_{L_+^2}^2\;dxdydt^\prime,
\end{split}
\end{equation}
and
\begin{equation}\label{416}
\begin{split}
-\int_0^t& \hbar(t')(e^\Psi\Delta_k^h\partial_y^2\tilde{u}_\Phi |e^\Psi\Delta_k^h\tilde{u}_\Phi)_{L_+^2}\;dt^\prime\\
&=\|\hbar^{\frac12}e^\Psi\Delta_k^h\partial_y\tilde{u}_\Phi(t)\|_{L^2(0, t;L_+^2)}^2+2\int_{0}^t\int_{\mathbb{R}_+^2} \hbar\partial_y\Psi e^{2\Psi}\Delta_k^h\tilde{u}_\Phi\Delta_k^h\partial_y\tilde{u}_\Phi\;dxdydt^\prime\\
&\ge\frac12\|\hbar^{\frac12}e^\Psi\Delta_k^h\partial_y\tilde{u}_\Phi(t)\|_{L^2(0, t;L_+^2)}^2
-2\int_{0}^t\int_{\mathbb{R}_+^2} \hbar(\partial_y\Psi)^2 |e^{\Psi}\Delta_k^h\tilde{u}_\Phi |^2\;dxdydt^\prime,
\end{split}
\end{equation}
in the last inequality we used Young's inequality.
Thanks to (\ref{214}), we obtain
\begin{equation}\label{417}
\begin{split}
\int_0^t &\hbar(t')(e^\Psi\Delta_k^h(\partial_t\tilde{u}_\Phi-\partial_y^2\tilde{u}_\Phi) |e^\Psi\Delta_k^h\tilde{u}_\Phi)_{L_+^2}\;dt^\prime\\
&\ge\frac12\left(\|\hbar^{\frac12}e^\Psi\Delta_k^h\tilde{u}_\Phi(t)\|_{L_+^2}^2
-\|\hbar^{\frac12}e^\Psi\Delta_k^h\tilde{u}_\Phi(0)\|_{L_+^2}^2\right.\\
&\left.\quad-\int_0^t \hbar^\prime(t^\prime)\|e^\Psi\Delta_k^h\tilde{u}_\Phi(t^\prime)\|_{L_+^2}^2\;dt^\prime
+\|\hbar^{\frac12}e^\Psi\Delta_k^h\partial_y\tilde{u}_\Phi(t)\|_{L^2(0, t;L_+^2)}^2\right).
\end{split}
\end{equation}
Similarly,
\begin{equation}\label{418}
\begin{split}
\int_0^t &\hbar(t')(e^\Psi\Delta_k^h(\partial_t\tilde{b}_\Phi-\partial_y^2\tilde{b}_\Phi) |e^\Psi\Delta_k^h\tilde{b}_\Phi)_{L_+^2}\;dt^\prime\\
&\ge\frac12\left(\|\hbar^{\frac12}e^\Psi\Delta_k^h\tilde{b}_\Phi(t)\|_{L_+^2}^2
-\|\hbar^{\frac12}e^\Psi\Delta_k^h\tilde{b}_\Phi(0)\|_{L_+^2}^2\right.\\
&\left.\quad-\int_0^t \hbar^\prime(t^\prime)\|e^\Psi\Delta_k^h\tilde{b}_\Phi(t^\prime)\|_{L_+^2}^2\;dt^\prime
+\|\hbar^{\frac12}e^\Psi\Delta_k^h\partial_y\tilde{b}_\Phi(t)\|_{L^2(0, t;L_+^2)}^2\right).
\end{split}
\end{equation}
$\bullet$ Estimate of $\int_0^t\lambda\dot{\eta}(t')\hbar(t')(e^\Psi |D_h|\Delta_k^h\tilde{u}_\Phi |e^\Psi\Delta_k^h\tilde{u}_\Phi)_{L_+^2}\;dt^\prime$ and

\qquad\qquad$\int_0^t\lambda\dot{\eta}(t')\hbar(t')(e^\Psi |D_h|\Delta_k^h\tilde{b}_\Phi |e^\Psi\Delta_k^h\tilde{b}_\Phi)_{L_+^2}\;dt^\prime$\\

From Lemma \ref{L2}, we immediately get that
\begin{equation}\label{419}
\int_0^t\lambda\dot{\eta}(t')\hbar(t')(e^\Psi |D_h|\Delta_k^h\tilde{u}_\Phi |e^\Psi\Delta_k^h\tilde{u}_\Phi)_{L_+^2}\;dt^\prime
\ge c\lambda 2^k\int_0^t\dot{\eta}(t')\|\hbar^{\frac12}e^\Psi\Delta_k^h\tilde{u}_\Phi\|_{L_+^2}^2\;dt^\prime,
\end{equation}
and
\begin{equation}\label{420}
\int_0^t\lambda\dot{\eta}(t')\hbar(t')(e^\Psi |D_h|\Delta_k^h\tilde{b}_\Phi |e^\Psi\Delta_k^h\tilde{b}_\Phi)_{L_+^2}\;dt^\prime
\ge c\lambda 2^k\int_0^t\dot{\eta}(t')\|\hbar^{\frac12}e^\Psi\Delta_k^h\tilde{b}_\Phi\|_{L_+^2}^2\;dt^\prime.
\end{equation}
$\bullet$ Estimate of $\int_0^t \hbar(t')(e^\Psi\Delta_k^h[(u+u^s+f(t)\varrho(y))\partial_x\tilde{u}]_\Phi |e^\Psi\Delta_k^h\tilde{u}_\Phi)_{L_+^2}\;dt^\prime$ and

\qquad\qquad $\int_0^t \hbar(t')(e^\Psi\Delta_k^h[(u+u^s+f(t)\varrho(y))\partial_x\tilde{b}]_\Phi |e^\Psi\Delta_k^h\tilde{b}_\Phi)_{L_+^2}\;dt^\prime$\\

Integration by parts implies that
\begin{equation}\label{421}
\int_0^t \hbar(t')(e^\Psi\Delta_k^h[(u^s+f(t)\varrho(y))\partial_x\tilde{u}]_\Phi |e^\Psi\Delta_k^h\tilde{u}_\Phi)_{L_+^2}\;dt^\prime=0.
\end{equation}
Then, applying Bony's decomposition to the remaining term $u\partial_x\tilde{u}$ gives
\[
u\partial_x\tilde{u}=T_u^h\partial_x\tilde{u}+T_{\partial_x\tilde{u}}^hu+R^h(u, \partial_x\tilde{u}).
\]
By virtue of (\ref{213}) and  the support properties of the Fourier transform of the terms in $T_u^h\partial_x\tilde{u}$, we have
\begin{align*}
&\int_0^t \hbar(t')|(e^\Psi\Delta_k^h[T_u^h\partial_x\tilde{u}]_\Phi |e^\Psi\Delta_k^h\tilde{u}_\Phi)_{L_+^2}|\;dt^\prime\\
\lesssim&\sum\limits_{|k'-k|\le4}\int_0^t \|S_{k'-1}^hu
_\Phi(t')\|_{L_+^\infty}\|\hbar^{\frac12}e^\Psi\Delta_{k'}^h\partial_x\tilde{u}_\Phi(t')\|_{L_+^2}\|\hbar^{\frac12}e^\Psi\Delta_k^h\tilde{u}_\Phi(t')\|_{L_+^2}\;dt'.
\end{align*}
Notice that
\begin{equation}\label{422}
\begin{split}
\|\Delta_k^hu_\Phi(t')\|_{L_v^\infty(L_h^2)}&\lesssim\|\Delta_k^h(\int_0^y\partial_yu_\Phi(t')\;dy')\|_{L_v^\infty(L_h^2)}\\
&\lesssim\|\Delta_k^h\partial_yu_\Phi(t')\|_{L_v^1(L_h^2)}\\
&\lesssim\|e^{-\Psi}\|_{L_v^2} \|e^\Psi\Delta_k^h\partial_yu_\Phi(t')\|_{L_+^2}\\
&\lesssim d_k(t')2^{-\frac k2}\langle t'\rangle^{\frac14}\|e^\Psi\partial_yu_\Phi(t')\|_{\mathcal{B}^{\frac12, 0}},
\end{split}
\end{equation}
where $\{d_k(t')\}_{k\in\mathbb{Z}}$ denotes a non-negative generic element in the unit sphere of $\ell^1(\mathbb{Z})$ for any $t'>0$ and satisfy $\sum\limits_{k\in\mathbb{Z}}d_k=1$. Then by Lemma \ref{L2}, we obtain
\begin{align*}
\|S_{k'-1}^hu_\Phi(t')\|_{L_+^\infty}&\lesssim\sum\limits_{k\le k'-2} 2^{\frac k2}\|\Delta_k^hu_\Phi(t')|_{L_v^\infty(L_h^2)}\\
&\lesssim\langle t'\rangle^{\frac14}\|e^\Psi\partial_yu_\Phi(t')\|_{\mathcal{B}^{\frac12, 0}}\lesssim\dot{\eta}(t').
\end{align*}
In this way, by applying H\"older's inequality, we get
\begin{equation}\label{423}
\begin{split}
&\int_0^t \hbar(t')|(e^\Psi\Delta_k^h[T_u^h\partial_x\tilde{u}]_\Phi |e^\Psi\Delta_k^h\tilde{u}_\Phi)_{L_+^2}|\;dt^\prime\\
\lesssim& 2^{k'}\left(\int_0^t\dot{\eta}(t')\|\hbar^{\frac12}e^\Psi\Delta_{k'}^h\tilde{u}_\Phi(t')\|_{L_+^2}^2\;dt'\right)^{\frac12}
\left(\int_0^t\dot{\eta}(t')\|\hbar^{\frac12}e^\Psi\Delta_{k}^h\tilde{u}_\Phi(t')\|_{L_+^2}^2\;dt'\right)^{\frac12}\\
\lesssim& d_k^22^{-k}\|\hbar^{\frac12}e^\Psi\tilde{u}_\Phi\|_{\tilde{L}_{t, \dot{\eta}(t)}^2(\mathcal{B}^{1,0})}^2.
\end{split}
\end{equation}
Similarly, by using (\ref{422}), Lemma \ref{L2} and H\"older's inequality again, we have
\begin{align*}
&\int_0^t \hbar(t)|(e^\Psi\Delta_k^h[T_{\partial_x\tilde{u}}^hu]_\Phi |e^\Psi\Delta_k^h\tilde{u}_\Phi)_{L_+^2}|\;dt^\prime\\
\lesssim&\sum\limits_{|k'-k|\le4}\int_0^t \|\hbar^{\frac12}e^\Psi S_{k'-1}^h\partial_x\tilde{u}
_\Phi(t')\|_{L_v^2(L_h^\infty)}
\|\Delta_{k'}^hu_\Phi(t')\|_{L_v^\infty(L_h^2)}
\|\hbar^{\frac12}e^\Psi\Delta_k^h\tilde{u}_\Phi(t')\|_{L_+^2}\;dt'\\
\lesssim&\sum\limits_{|k'-k|\le4} d_{k'}2^{-\frac {k'}2}\left(\int_0^t\dot{\eta}(t')\|\hbar^{\frac12}e^\Psi S_{k'-1}^h\partial_x\tilde{u}_\Phi(t')\|_{L_v^2(L_h^\infty)}^2\;dt'\right)^{\frac12}\left(\int_0^t\dot{\eta}(t')\|\hbar^{\frac12}e^\Psi\Delta_{k}^h\tilde{u}_\Phi(t')\|_{L_+^2}^2\;dt'\right)^{\frac12}.
\end{align*}
By Lemma \ref{L2}, it yields
\begin{align*}
\left(\int_0^t\dot{\eta}(t')\|\hbar^{\frac12}e^\Psi S_{k'-1}^h\partial_x\tilde{u}_\Phi(t')\|_{L_v^2(L_h^\infty)}^2\;dt'\right)^{\frac12}
&\lesssim\sum\limits_{j\le k'-2}2^{\frac{3j}2}\left(\int_0^t\dot{\eta}(t')\|\hbar^{\frac12}e^\Psi\Delta_{j}^h\tilde{u}_\Phi(t')\|_{L_+^2}^2\;dt'\right)^{\frac12}\\
&\lesssim 2^{\frac{k'}2}\|\hbar^{\frac12}e^\Psi\tilde{u}_\Phi\|_{\tilde{L}_{t, \dot{\eta}(t)}^2(\mathcal{B}^{1,0})}.
\end{align*}
Hence, we deduce that
\begin{equation}\label{424}
\int_0^t \hbar(t')|(e^\Psi\Delta_k^h[T_{\partial_x\tilde{u}}^hu]_\Phi |e^\Psi\Delta_k^h\tilde{u}_\Phi)_{L_+^2}|\;dt^\prime
\lesssim d_k^22^{-k}\|\hbar^{\frac12}e^\Psi\tilde{u}_\Phi\|_{\tilde{L}_{t, \dot{\eta}(t)}^2(\mathcal{B}^{1,0})}^2.
\end{equation}
For the last term, by Lemma \ref{L2} and (\ref{422}) again, it holds that
{\begin{align*}
&\int_0^t \hbar(t')|(e^\Psi\Delta_k^h[R^h(u, \partial_x\tilde{u})]_\Phi |e^\Psi\Delta_k^h\tilde{u}_\Phi)_{L_+^2}|\;dt^\prime\\
\lesssim&\sum\limits_{k'\ge k-3}\int_0^t \|\tilde{\Delta}_{k'}^hu_\Phi(t')\|_{L_+^\infty}
\|\hbar^{\frac12}e^\Psi\Delta_{k'}^h\partial_x\tilde{u}_\Phi(t')\|_{L_+^2}
\|\hbar^{\frac12}e^\Psi\Delta_k^h\tilde{u}_\Phi(t')\|_{L_+^2}\;dt'\\
\lesssim&\sum\limits_{k'\ge k-3}2^{k'}\int_0^t \dot{\eta}(t')
\|\hbar^{\frac12}e^\Psi\Delta_{k'}^h\tilde{u}_\Phi(t')\|_{L_+^2}
\|\hbar^{\frac12}e^\Psi\Delta_k^h\tilde{u}_\Phi(t')\|_{L_+^2}\;dt',
\end{align*}
then
\begin{equation}\label{425}
\int_0^t \hbar(t')|(e^\Psi\Delta_k^h[R^h(u, \partial_x\tilde{u})]_\Phi |e^\Psi\Delta_k^h\tilde{u}_\Phi)_{L_+^2}|\;dt^\prime\lesssim d_k^22^{-k}\|\hbar^{\frac12}e^\Psi\tilde{u}_\Phi\|_{\tilde{L}_{t, \dot{\eta}(t)}^2(\mathcal{B}^{1,0})}^2.
\end{equation}}
Combining (\ref{421}) and (\ref{423})-(\ref{425}) together, we get the estimate
\begin{equation}\label{426}
\int_0^t \hbar(t')|(e^\Psi\Delta_k^h[(u+u^s+f(t)\varrho(y))\partial_x\tilde{u}]_\Phi |e^\Psi\Delta_k^h\tilde{u}_\Phi)_{L_+^2}|\;dt^\prime\lesssim d_k^22^{-k}\|\hbar^{\frac12}e^\Psi\tilde{u}_\Phi\|_{\tilde{L}_{t, \dot{\eta}(t)}^2(\mathcal{B}^{1,0})}^2.
\end{equation}
By similar arguments, we have
\begin{equation}\label{427}
\int_0^t \hbar(t')|(e^\Psi\Delta_k^h[(u+u^s+f(t)\varrho(y))\partial_x\tilde{b}]_\Phi |e^\Psi\Delta_k^h\tilde{b}_\Phi)_{L_+^2}|\;dt^\prime\lesssim d_k^22^{-k}\|\hbar^{\frac12}e^\Psi\tilde{b}_\Phi\|_{\tilde{L}_{t, \dot{\eta}(t)}^2(\mathcal{B}^{1,0})}^2.
\end{equation}
$\bullet$ Estimate of $\int_0^t\hbar(t')(e^\Psi \Delta_k^h[(1+b)\partial_x\tilde{b}]_\Phi |e^\Psi\Delta_k^h\tilde{u}_\Phi)_{L_+^2}\;dt^\prime$ and

 \qquad\qquad $\int_0^t\hbar(t')(e^\Psi \Delta_k^h[(1+b)\partial_x\tilde{u}]_\Phi |e^\Psi\Delta_k^h\tilde{b}_\Phi)_{L_+^2}\;dt^\prime$\\

Indeed,  the estimates of these two terms are similar to the estimates above. We only give the estimates and omit the details here for simplicity.
\begin{equation}\label{428}
\begin{split}
&\int_0^t \hbar(t')|(e^\Psi\Delta_k^h[(1+b)\partial_x\tilde{b}]_\Phi |e^\Psi\Delta_k^h\tilde{u}_\Phi)_{L_+^2}|\;dt^\prime\\
\lesssim& d_k^22^{-k}\left(\|\hbar^{\frac12}e^\Psi\tilde{u}_\Phi\|_{\tilde{L}_{t, \dot{\eta}(t)}^2(\mathcal{B}^{1,0})}^2+\|\hbar^{\frac12}e^\Psi\tilde{b}_\Phi\|_{\tilde{L}_{t, \dot{\eta}(t)}^2(\mathcal{B}^{1,0})}^2\right),
\end{split}
\end{equation}
and
\begin{equation}\label{429}
\begin{split}
&\int_0^t \hbar(t')|(e^\Psi\Delta_k^h[(1+b)\partial_x\tilde{u}]_\Phi |e^\Psi\Delta_k^h\tilde{b}_\Phi)_{L_+^2}|\;dt^\prime\\
\lesssim& d_k^22^{-k}\left(\|\hbar^{\frac12}e^\Psi\tilde{u}_\Phi\|_{\tilde{L}_{t, \dot{\eta}(t)}^2(\mathcal{B}^{1,0})}^2+\|\hbar^{\frac12}e^\Psi\tilde{b}_\Phi\|_{\tilde{L}_{t, \dot{\eta}(t)}^2(\mathcal{B}^{1,0})}^2\right).
\end{split}
\end{equation}
$\bullet$ Estimate of $\int_0^t\hbar(t')(e^\Psi \Delta_k^h[v\partial_y\tilde{u}]_\Phi |e^\Psi\Delta_k^h\tilde{u}_\Phi)_{L_+^2}\;dt^\prime$, $\int_0^t\hbar(t')(e^\Psi \Delta_k^h[v\partial_y\tilde{b}]_\Phi |e^\Psi\Delta_k^h\tilde{b}_\Phi)_{L_+^2}\;dt^\prime$ and

\qquad\qquad$\int_0^t\hbar(t')(e^\Psi \Delta_k^h[g\partial_y\tilde{b}]_\Phi |e^\Psi\Delta_k^h\tilde{u}_\Phi)_{L_+^2}\;dt^\prime$, $\int_0^t\hbar(t')(e^\Psi \Delta_k^h[g\partial_y\tilde{u}]_\Phi |e^\Psi\Delta_k^h\tilde{b}_\Phi)_{L_+^2}\;dt^\prime$

Again, by Bony's decomposition for $v\partial_y\tilde{u}$, it gives
\[
v\partial_y\tilde{u}=T_v^h\partial_y\tilde{u}+T_{\partial_y\tilde{u}}^hv+R^h(v, \partial_y\tilde{u}).
\]
It derives from a direct calculation that
\begin{align*}
&\int_0^t \hbar(t')|(e^\Psi\Delta_k^h[T_v^h\partial_y\tilde{u}]_\Phi |e^\Psi\Delta_k^h\tilde{u}_\Phi)_{L_+^2}|\;dt^\prime\\
\lesssim&\sum\limits_{|k'-k|\le4}\int_0^t \|\hbar^{\frac12}S_{k'-1}^hv_\Phi(t')\|_{L_+^\infty}
\|e^\Psi\Delta_{k'}^h\partial_y\tilde{u}_\Phi(t')\|_{L_+^2}\|\hbar^{\frac12}e^\Psi\Delta_k^h\tilde{u}_\Phi(t')\|_{L_+^2}\;dt'\\
\lesssim&\sum\limits_{|k'-k|\le4}d_{k'} 2^{-\frac{k'}2}\int_0^t \|\hbar^{\frac12}S_{k'-1}^h(\int_0^y\partial_xu_\Phi(t')\;dy')\|_{L_+^\infty}
\|e^\Psi\partial_y\tilde{u}_\Phi(t')\|_{\mathcal{B}^{\frac12, 0}}\|\hbar^{\frac12}e^\Psi\Delta_k^h\tilde{u}_\Phi(t')\|_{L_+^2}\;dt'.
\end{align*}
On the other hand, by the similar derivation of (\ref{422}), we have
\begin{align*}
\|S_{k'-1}^h(\int_0^y\partial_xu_\Phi(t')\;dy')\|_{L_+^\infty}&\lesssim\sum\limits_{k\le k'-2}2^{\frac32 k}\|\Delta_{k}^hu_\Phi (t')\|_{L_v^1(L_h^2)}\\
&\lesssim\sum\limits_{k\le k'-2}2^{\frac32 k}\|e^{-\Psi} (t')\|_{L_v^2}\|e^\Psi\Delta_{k}^hu_\Phi (t')\|_{L_+^2}\\
&\lesssim\sum\limits_{k\le k'-2}2^{\frac32 k}\langle t'\rangle^{\frac14}\|e^\Psi\Delta_{k}^hu_\Phi (t')\|_{L_+^2},
\end{align*}
from which and the above inequality, we arrive at
\begin{equation}\label{430}
\int_0^t \hbar(t')|(e^\Psi\Delta_k^h[T_v^h\partial_y\tilde{u}]_\Phi |e^\Psi\Delta_k^h\tilde{u}_\Phi)_{L_+^2}|\;dt^\prime\lesssim d_k^22^{-k}\|\hbar^{\frac12}e^\Psi\tilde{u}_\Phi\|_{\tilde{L}_{t, \dot{\eta}(t)}^2(\mathcal{B}^{1,0})}^2.
\end{equation}
By the same manner, we have
\begin{align*}
&\int_0^t \hbar(t')|(e^\Psi\Delta_k^h[T_{\partial_y\tilde{u}}^hv]_\Phi |e^\Psi\Delta_k^h\tilde{u}_\Phi)_{L_+^2}|\;dt^\prime\\
\lesssim&\sum\limits_{|k'-k|\le4}\int_0^t
\|e^\Psi S_{k'-1}^h\partial_y\tilde{u}_\Phi(t')\|_{L_v^2(L_h^\infty)}
\|\hbar^{\frac12}\Delta_{k'}^h(\int_0^y \partial_yv_\Phi(t')\;dy')\|_{L_v^\infty(L_h^2)}
\|\hbar^{\frac12}e^\Psi\Delta_k^h\tilde{u}_\Phi(t')\|_{L_+^2}\;dt'.
\end{align*}
Again, we deduce by a similar derivation of (\ref{422}) that
\begin{equation*}
\|e^\Psi S_{k'-1}^h\partial_y\tilde{u}_\Phi(t')\|_{L_v^2(L_h^\infty)}\lesssim\sum\limits_{k\le k'-2}2^{\frac{k}2}\|e^\Psi\Delta_{k}^h\partial_y\tilde{u}_\Phi (t')\|_{L_+^2}\lesssim\|e^{\Psi} \partial_y\tilde{u}_\Phi (t')\|_{\mathcal{B}^{\frac12, 0}}
\end{equation*}
and
\begin{equation}\label{431}
\begin{split}
\|\hbar^{\frac12}\Delta_{k'}^h(\int_0^y \partial_yv_\Phi(t')\;dy')\|_{L_v^\infty(L_h^2)}&\lesssim
2^{k'}\|\hbar^{\frac12}\Delta_{k'}^h(\int_0^y u_\Phi(t')\;dy')\|_{L_v^\infty(L_h^2)}\\
&\lesssim
2^{k'}\|\hbar^{\frac12}\Delta_{k'}^h u_\Phi(t')\|_{L_v^1(L_h^2)}\\
&\lesssim
2^{k'}\|e^{-\Psi}(t')\|_{L_v^2}\|\hbar^{\frac12}e^\Psi\Delta_{k'}^h u_\Phi(t')\|_{L_+^2}\\
&\lesssim 2^{k'}\langle t'\rangle ^{\frac14}\|\hbar^{\frac12}e^\Psi\Delta_{k'}^h u_\Phi(t')\|_{L_+^2}.
\end{split}
\end{equation}
Consequently, by collecting the above three inequalities together, we obtain
\begin{equation}\label{432}
\int_0^t \hbar(t')|(e^\Psi\Delta_k^h[T_{\partial_y\tilde{u}}^hv]_\Phi |e^\Psi\Delta_k^h\tilde{u}_\Phi)_{L_+^2}|\;dt^\prime\lesssim d_k^22^{-k}\|\hbar^{\frac12}e^\Psi\tilde{u}_\Phi\|_{\tilde{L}_{t, \dot{\eta}(t)}^2(\mathcal{B}^{1,0})}^2.
\end{equation}
For the remaining term, it is estimated as follows.
\begin{align*}
&\int_0^t \hbar(t')|(e^\Psi\Delta_k^h[R^h(v, \partial_y\tilde{u})]_\Phi |e^\Psi\Delta_k^h\tilde{u}_\Phi)_{L_+^2}|\;dt^\prime\\
\lesssim&\sum\limits_{k'\ge k-3}\int_0^t
\|\hbar^{\frac12}\Delta_{k'}^h(\int_0^y\partial_yv_\Phi(t')\;dy')\|_{L_v^\infty(L_h^2)}
\|e^\Psi\tilde{\Delta}_{k'}^h \partial_y\tilde{u}_\Phi(t')\|_{L_v^2(L_h^\infty)}
\|\hbar^{\frac12}e^\Psi\Delta_k^h\tilde{u}_\Phi(t')\|_{L_+^2}\;dt'.
\end{align*}
From (\ref{431}) and the fact that
\[
\|e^\Psi\tilde{\Delta}_{k'}^h\partial_y\tilde{u}_\Phi(t')\|_{L_v^2(L_h^\infty)}\lesssim\sum\limits_{|k-k'|\le1} 2^{\frac k2}\|e^\Psi\Delta_{k}^h\partial_y\tilde{u}_\Phi(t')\|_{L_+^2},
\]
we get
\begin{equation}\label{433}
\int_0^t \hbar(t')|(e^\Psi\Delta_k^h[R^h(v, \partial_y\tilde{u})]_\Phi |e^\Psi\Delta_k^h\tilde{u}_\Phi)_{L_+^2}|\;dt^\prime\lesssim d_k^22^{-k}\|\hbar^{\frac12}e^\Psi\tilde{u}_\Phi\|_{\tilde{L}_{t, \dot{\eta}(t)}^2(\mathcal{B}^{1,0})}^2.
\end{equation}
As a sequence, the following estimate holds true.
\begin{equation}\label{434}
\int_0^t \hbar(t')|(e^\Psi\Delta_k^h[v\partial_y\tilde{u}]_\Phi |e^\Psi\Delta_k^h\tilde{u}_\Phi)_{L_+^2}|\;dt^\prime\lesssim d_k^22^{-k}\|\hbar^{\frac12}e^\Psi\tilde{u}_\Phi\|_{\tilde{L}_{t, \dot{\eta}(t)}^2(\mathcal{B}^{1,0})}^2.
\end{equation}
Along the same line, we also have
\begin{equation}\label{435}
\begin{split}
\int_0^t\hbar(t')|(e^\Psi &\Delta_k^h[v\partial_y\tilde{b}]_\Phi |e^\Psi\Delta_k^h\tilde{b}_\Phi)_{L_+^2}|\;dt^\prime\\
&\lesssim d_k^22^{-k}\left(\|\hbar^{\frac12}e^\Psi\tilde{u}_\Phi\|_{\tilde{L}_{t, \dot{\eta}(t)}^2(\mathcal{B}^{1,0})}^2+\|\hbar^{\frac12}e^\Psi\tilde{b}_\Phi\|_{\tilde{L}_{t, \dot{\eta}(t)}^2(\mathcal{B}^{1,0})}^2\right),
\end{split}
\end{equation}
\begin{equation}\label{436}
\begin{split}
\int_0^t\hbar(t')|(e^\Psi &\Delta_k^h[g\partial_y\tilde{b}]_\Phi |e^\Psi\Delta_k^h\tilde{u}_\Phi)_{L_+^2}|\;dt^\prime\\
&\lesssim d_k^22^{-k}\left(\|\hbar^{\frac12}e^\Psi\tilde{u}_\Phi\|_{\tilde{L}_{t, \dot{\eta}(t)}^2(\mathcal{B}^{1,0})}^2+\|\hbar^{\frac12}e^\Psi\tilde{b}_\Phi\|_{\tilde{L}_{t, \dot{\eta}(t)}^2(\mathcal{B}^{1,0})}^2\right),
\end{split}
\end{equation}
and
\begin{equation}\label{437}
\int_0^t \hbar(t')|(e^\Psi\Delta_k^h[g\partial_y\tilde{b}]_\Phi |e^\Psi\Delta_k^h\tilde{b}_\Phi)_{L_+^2}|\;dt^\prime\lesssim d_k^22^{-k}\|\hbar^{\frac12}e^\Psi\tilde{b}_\Phi\|_{\tilde{L}_{t, \dot{\eta}(t)}^2(\mathcal{B}^{1,0})}^2.
\end{equation}
$\bullet$ Estimate of $\int_0^t \hbar(t')(\partial_y^2u^se^\Psi \Delta_k^h[v\psi]_\Phi |e^\Psi\Delta_k^h\tilde{u}_\Phi)_{L_+^2}\;dt^\prime$ and $\int_0^t \hbar(t')(\partial_y^2u^se^\Psi \Delta_k^h[g\psi]_\Phi |e^\Psi\Delta_k^h\tilde{b}_\Phi)_{L_+^2}\;dt^\prime$\\

Applying Bony's decomposition to $v\psi$ to get
\[
v\psi=T_v^h\psi+T_\psi^hv+R^h(v, \psi).
\]
Considering (\ref{213}) and the support properties of the Fourier transform of the terms in $T_v^h\psi$, we have
\begin{align*}
&\int_0^t \hbar(t')|(\partial_y^2u^se^\Psi \Delta_k^h[T_v^h\psi]_\Phi |e^\Psi\Delta_k^h\tilde{u}_\Phi)_{L_+^2}|\;dt^\prime\\
\lesssim&\sum\limits_{|k'-k|\le 4}\int_0^t
\|S_{k'-1}^hv_\Phi(t')\|_{L_v^\infty(L_h^2)}
\|\hbar^{\frac12}\Delta_{k'}^h\psi_\Phi(t')\|_{L_+^\infty}
\|e^\Psi\partial_y^2u^s\|_{L_v^2}
\|\hbar^{\frac12}e^\Psi\Delta_k^h\tilde{u}_\Phi(t')\|_{L_+^2}\;dt'.
\end{align*}
Notice that the following estimates are direct consequences of Lemmas \ref{L2} and \ref{L2.1}.
\begin{align*}
\|S_{k'-1}^hv_\Phi(t')\|_{L_v^\infty(L_h^2)}
&\lesssim\sum\limits_{k\le k'-2}
\|\Delta_{k}^h (\int_0^y\partial_yv_\Phi(t')\;dy')\|_{L_v^\infty(L_h^2)}\\
&\lesssim\sum\limits_{k\le k'-2}2^{k}
\|\Delta_{k}^h u_\Phi(t')\|_{L_v^1(L_h^2)}\\
&\lesssim\sum\limits_{k\le k'-2}2^{k}
\|e^{-\Psi}(t')\|_{L_v^2}\|e^\Psi\Delta_{k}^h u_\Phi(t')\|_{L_+^2}\\
&\lesssim\sum\limits_{k\le k'-2}
 2^{k} \langle t'\rangle ^{\frac14}\|e^\Psi\Delta_{k}^h u_\Phi(t')\|_{L_+^2}\\
 &\lesssim\sum\limits_{k\le k'-2}
 2^{k} \langle t'\rangle ^{\frac34}\|e^\Psi\Delta_{k}^h \partial_yu_\Phi(t')\|_{L_+^2}\\
 &\lesssim
 \langle t'\rangle ^{\frac12}\dot{\eta}(t'),
 \end{align*}
and
\begin{equation}\label{438}
\begin{split}
\|\hbar^{\frac12}\Delta_{k'}^h\psi_\Phi(t')\|_{L_+^\infty}
&\lesssim
\|\hbar^{\frac12}\Delta_{k'}^h (\int_0^y\partial_y\psi_\Phi(t')\;dy')\|_{L_v^\infty(L_h^2)}\\
&\lesssim
\|\hbar^{\frac12}\Delta_{k'}^h \tilde{b}_\Phi(t')\|_{L_v^1(L_h^2)}\\
&\lesssim
\|e^{-\Psi}(t')\|_{L_v^2}\|\hbar^{\frac12}e^\Psi\Delta_{k'}^h \tilde{b}_\Phi(t')\|_{L_+^2}\\
&\lesssim \langle t'\rangle ^{\frac14}\|\hbar^{\frac12}e^\Psi\Delta_{k'}^h \tilde{b}_\Phi(t')\|_{L_+^2}.
\end{split}
\end{equation}
Then, by using $\|e^\Psi\partial_y^2u^s(t, \cdot)\|_{L_v^2}\le \frac C{{\langle t\rangle}^{3/4}}$ in (\ref{116}), we obtain
\begin{equation}\label{439}
\begin{split}
&\int_0^t \hbar(t')|(\partial_y^2u^se^\Psi \Delta_k^h[T_v^h\psi]_\Phi |e^\Psi\Delta_k^h\tilde{u}_\Phi)_{L_+^2}|\;dt^\prime\\
\lesssim& d_k^22^{-k}\left(\|\hbar^{\frac12}e^\Psi\tilde{u}_\Phi\|_{\tilde{L}_{t, \dot{\eta}(t)}^2(\mathcal{B}^{1,0})}^2+\|\hbar^{\frac12}e^\Psi\tilde{b}_\Phi\|_{\tilde{L}_{t, \dot{\eta}(t)}^2(\mathcal{B}^{1,0})}^2\right).
\end{split}
\end{equation}
Similarly,
\begin{align*}
&\int_0^t \hbar(t')|(\partial_y^2u^se^\Psi \Delta_k^h[T_\psi^hv]_\Phi |e^\Psi\Delta_k^h\tilde{u}_\Phi)_{L_+^2}|\;dt^\prime\\
\lesssim&\sum\limits_{|k'-k|\le 4}\int_0^t
\|\hbar^{\frac12}S_{k'-1}^h\psi_\Phi(t')\|_{L_v^\infty(L_h^2)}
\|\Delta_{k'}^hv_\Phi(t')\|_{L_+^\infty}
\|e^\Psi\partial_y^2u^s\|_{L_v^2}
\|\hbar^{\frac12}e^\Psi\Delta_k^h\tilde{u}_\Phi(t')\|_{L_+^2}\;dt',
\end{align*}
where
\begin{align*}
\|\hbar^{\frac12}S_{k'-1}^h\psi_\Phi(t')\|_{L_v^\infty(L_h^2)}&\lesssim\sum\limits_{k\le k'-2}
\|\hbar^{\frac12}\Delta_{k}^h\tilde{b}(t')\|_{L_v^1(L_h^2)}\\
&\lesssim\sum\limits_{k\le k'-2}
\langle t'\rangle ^{\frac14}\|\hbar^{\frac12}e^\Psi\Delta_{k}^h\tilde{b}(t')
\|_{L_+^2}
\end{align*}
and
\begin{align*}
\|\Delta_{k'}^hv_\Phi(t')\|_{L_+^\infty}
&\lesssim
\|\Delta_{k'}^h (\int_0^y\partial_xu_\Phi(t')\;dy')\|_{L_+^\infty}\\
&\lesssim
2^{\frac{3k'}2}\|\Delta_{k'}^hu_\Phi(t')\|_{L_v^1(L_h^2)}\\
&\lesssim
2^{\frac{3k'}2}\|e^{-\Psi}(t')\|_{L_v^2}\|e^\Psi\Delta_{k'}^h u_\Phi(t')\|_{L_+^2}\\
&\lesssim
 2^{\frac{3k'}2} \langle t'\rangle ^{\frac34}\|e^\Psi\Delta_{k'}^h \partial_yu_\Phi(t')\|_{L_+^2}\\
 &\lesssim
 d_{k'}(t')2^{k'} \langle t'\rangle ^{\frac34}\|e^\Psi\partial_yu_\Phi(t')\|_{\mathcal{B}^{\frac12, 0}}.
 \end{align*}
Finally, we get the following estimate.
\begin{equation}\label{440}
\begin{split}
&\int_0^t \hbar(t')|(\partial_y^2u^se^\Psi \Delta_k^h[T_\psi^hv]_\Phi |e^\Psi\Delta_k^h\tilde{u}_\Phi)_{L_+^2}|\;dt^\prime\\
\lesssim& d_k^22^{-k}\left(\|\hbar^{\frac12}e^\Psi\tilde{u}_\Phi\|_{\tilde{L}_{t, \dot{\eta}(t)}^2(\mathcal{B}^{1,0})}^2+\|\hbar^{\frac12}e^\Psi\tilde{b}_\Phi\|_{\tilde{L}_{t, \dot{\eta}(t)}^2(\mathcal{B}^{1,0})}^2\right).
\end{split}
\end{equation}
Using the support properties of the Fourier transform of the terms in $R^h(v, \psi)$, we obtain
\begin{align*}
&\int_0^t \hbar(t')|(\partial_y^2u^se^\Psi \Delta_k^h[R^h(v, \psi)]_\Phi |e^\Psi\Delta_k^h\tilde{u}_\Phi)_{L_+^2}|\;dt^\prime\\
\lesssim&\sum\limits_{k'\ge k-3}\int_0^t
\|\tilde{\Delta}_{k'}^hv_\Phi(t')\|_{L_v^\infty(L_h^2)}
\|\hbar^{\frac12}\Delta_{k}^h\psi_\Phi(t')\|_{L_+^\infty}
\|e^\Psi\partial_y^2u^s\|_{L_v^2}
\|\hbar^{\frac12}e^\Psi\Delta_k^h\tilde{u}_\Phi(t')\|_{L_+^2}\;dt'.
\end{align*}
By the similar arguments as in (\ref{431}), we infer
\begin{align*}
\|\tilde{\Delta}_{k'}^hv_\Phi(t')\|_{L_v^\infty(L_h^2)}&\lesssim\sum\limits_{|k-k'|\le 1}
2^{k}\|\Delta_{k}^hu_\Phi(t')\|_{L_v^1(L_h^2)}\\
&\lesssim\sum\limits_{|k-k'|\le 1}
2^{k}\langle t'\rangle^{\frac34}\|e^\Psi\Delta_{k}^h\partial_yu_\Phi(t')\|_{L_+^2}.
\end{align*}
From which and (\ref{438}), we have
\begin{equation}\label{441}
\begin{split}
&\int_0^t \hbar(t')|(\partial_y^2u^se^\Psi \Delta_k^h[R^h(v, \psi)]_\Phi |e^\Psi\Delta_k^h\tilde{u}_\Phi)_{L_+^2}|\;dt^\prime\\
\lesssim& d_k^22^{-k}\left(\|\hbar^{\frac12}e^\Psi\tilde{u}_\Phi\|_{\tilde{L}_{t, \dot{\eta}(t)}^2(\mathcal{B}^{1,0})}^2+\|\hbar^{\frac12}e^\Psi\tilde{b}_\Phi\|_{\tilde{L}_{t, \dot{\eta}(t)}^2(\mathcal{B}^{1,0})}^2\right).
\end{split}
\end{equation}
Hence, it follows, by summing up (\ref{439})-(\ref{441}), that
\begin{equation}\label{442}
\begin{split}
&\int_0^t \hbar(t')|(\partial_y^2u^se^\Psi \Delta_k^h[v \psi]_\Phi |e^\Psi\Delta_k^h\tilde{u}_\Phi)_{L_+^2}|\;dt^\prime\\
\lesssim& d_k^22^{-k}\left(\|\hbar^{\frac12}e^\Psi\tilde{u}_\Phi\|_{\tilde{L}_{t, \dot{\eta}(t)}^2(\mathcal{B}^{1,0})}^2+\|\hbar^{\frac12}e^\Psi\tilde{b}_\Phi\|_{\tilde{L}_{t, \dot{\eta}(t)}^2(\mathcal{B}^{1,0})}^2\right).
\end{split}
\end{equation}
In a similar way, it also holds
\begin{equation}\label{443}
\int_0^t \hbar(t')|(\partial_y^2u^se^\Psi \Delta_k^h[g \psi]_\Phi |e^\Psi\Delta_k^h\tilde{b}_\Phi)_{L_+^2}|\;dt^\prime
\lesssim d_k^22^{-k}\|\hbar^{\frac12}e^\Psi\tilde{b}_\Phi\|_{\tilde{L}_{t, \dot{\eta}(t)}^2(\mathcal{B}^{1,0})}^2.
\end{equation}
$\bullet$ Estimate of $\int_0^t \hbar(t')(f(t')\varrho'' e^\Psi \Delta_k^h[v\psi]_\Phi |e^\Psi\Delta_k^h\tilde{u}_\Phi)_{L_+^2}\;dt^\prime$ and

\qquad\qquad $\int_0^t \hbar(t')(f(t')\varrho'' e^\Psi \Delta_k^h[g\psi]_\Phi |e^\Psi\Delta_k^h\tilde{b}_\Phi)_{L_+^2}\;dt^\prime$\\

As a matter of fact in (\ref{115}),
\[
|f(t')|\le C\langle t'\rangle^{-1}\le C\langle t'\rangle^{-3/4}.
\]
Moreover, since $\varrho''(y)$ is supported in $[1, 2]$, we have
\[
\|f(t')e^\Psi\varrho''\|_{L_v^2}\le C\langle t'\rangle^{-3/4}.
\]
Then, by the similar arguments as above, we obtain
\begin{equation}\label{444}
\begin{split}
\int_0^t \hbar(t')&|(f(t')\varrho'' e^\Psi \Delta_k^h[v \psi]_\Phi |e^\Psi\Delta_k^h\tilde{u}_\Phi)_{L_+^2}|\;dt^\prime\\
&\lesssim d_k^22^{-k}\left(\|\hbar^{\frac12}e^\Psi\tilde{u}_\Phi\|_{\tilde{L}_{t, \dot{\eta}(t)}^2(\mathcal{B}^{1,0})}^2+\|\hbar^{\frac12}e^\Psi\tilde{b}_\Phi\|_{\tilde{L}_{t, \dot{\eta}(t)}^2(\mathcal{B}^{1,0})}^2\right).
\end{split}
\end{equation}
Likewise, we also have
\begin{equation}\label{445}
\int_0^t \hbar(t')|(f(t')\varrho''e^\Psi \Delta_k^h[g \psi]_\Phi |e^\Psi\Delta_k^h\tilde{b}_\Phi)_{L_+^2}|\;dt^\prime\lesssim d_k^22^{-k}\|\hbar^{\frac12}e^\Psi\tilde{b}_\Phi\|_{\tilde{L}_{t, \dot{\eta}(t)}^2(\mathcal{B}^{1,0})}^2.
\end{equation}
$\bullet$ Estimate of $\int_0^t \hbar(t')(2\partial_y^2u^se^\Psi \Delta_k^h\tilde{b}_\Phi |e^\Psi\Delta_k^h\tilde{u}_\Phi)_{L_+^2}\;dt^\prime$\\

By using $\|\partial_y^2u^s(t, \cdot)\|_{L_v^\infty}\le\frac C{\langle t\rangle}$ in (\ref{116}) and H\"older's inequality, we have for any $\alpha>0$
\begin{align*}
\int_0^t \hbar(t')|(2\partial_y^2u^se^\Psi \Delta_k^h\tilde{b}_\Phi |e^\Psi\Delta_k^h\tilde{u}_\Phi)_{L_+^2}|dt^\prime&\le\int_0^t 2C\hbar(t')\langle t\rangle^{-1}| (e^\Psi\Delta_k^h\tilde{b}_\Phi | e^\Psi\Delta_k^h\tilde{u}_\Phi)_{L_+^2}|\;dt'\\
&\le \int_0^t\langle t\rangle^{-1}\left(\alpha \|\hbar^{\frac12}e^\Psi\Delta_k^h\tilde{u}_\Phi\|_{L_+^2}^2+\frac{C^2}\alpha \|\hbar^{\frac12}e^\Psi\Delta_k^h\tilde{b}_\Phi\|_{L_+^2}^2\right)\;dt',
\end{align*}
which and Lemma \ref{L2.1} imply
\begin{align}\label{446}
&\int_0^t \hbar(t')|(2\partial_y^2u^se^\Psi \Delta_k^h\tilde{b}_\Phi |e^\Psi\Delta_k^h\tilde{u}_\Phi)_{L_+^2}|dt^\prime\nonumber\\
\le&
2\alpha \|\hbar^{\frac12}e^\Psi\Delta_k^h\partial_y\tilde{u}_\Phi\|_{L^2(0, t; L_+^2)}^2+\frac{2C^2}\alpha \|\hbar^{\frac12}e^\Psi\Delta_k^h\partial_y\tilde{b}_\Phi\|_{L^2(0, t; L_+^2)}^2.
\end{align}
$\bullet$Estimate of $\int_0^t \hbar(t')(2f(t')\varrho ''(y)e^\Psi \Delta_k^h\tilde{b}_\Phi |e^\Psi\Delta_k^h\tilde{u}_\Phi)_{L_+^2}\;dt^\prime$\\

Here, from (\ref{115}) and the fact that $\varrho''(y)$ is supported on the interval $[1, 2]$, we also have the estimate $\|f(t)\varrho'' (y)\|_{L_v^\infty}\le\frac C{\langle t\rangle}$. From which and Lemma \ref{L2.1}, one also has
\begin{align}
\label{447}
&\int_0^t \hbar(t')|(f(t')\varrho''e^\Psi \Delta_k^h\tilde{b}_\Phi |e^\Psi\Delta_k^h\tilde{u}_\Phi)_{L_+^2}|dt^\prime\nonumber\\
\le&
2\alpha \|\hbar^{\frac12}e^\Psi\Delta_k^h\partial_y\tilde{u}_\Phi\|_{L^2(0, t; L_+^2)}^2+\frac{2C^2}\alpha \|\hbar^{\frac12}e^\Psi\Delta_k^h\partial_y\tilde{b}_\Phi\|_{L^2(0, t; L_+^2)}^2.
\end{align}
Now, integrating both (\ref{413}) and (\ref{414}) over $[0, t]$ and inserting the above estimates into the resulting equalities, we conclude that
\begin{equation}\label{448}
\begin{split}
\|\hbar^{\frac12}&e^\Psi\Delta_k^h\tilde{u}_\Phi\|_{L^\infty(0, t; L_+^2)}^2
+2c\lambda 2^k\int_0^t \dot{\eta}(t')\|\hbar^{\frac12}e^\Psi\Delta_k^h\tilde{u}_\Phi(t')\|_{L_+^2}^2\;dt'
+\|\hbar^{\frac12}e^\Psi\Delta_k^h\partial_y\tilde{u}_\Phi\|_{L^2(0, t; L_+^2)}^2\\
&\le\|\hbar^{\frac12}e^\Psi\Delta_k^h\tilde{u}_\Phi(0)\|_{L_+^2}^2
+\int_0^t \hbar '(t')\|e^\Psi\Delta_k^h\tilde{u}_\Phi\|_{L_+^2}^2\;dt'\\
&\qquad+Cd_k^2 2^{-k}\left(\|\hbar^{\frac12}e^\Psi\tilde{u}_\Phi\|_{\tilde{L}^2_{t, \dot{\eta}(t)}(\mathcal{B}^{1, 0})}^2+\|\hbar^{\frac12}e^\Psi\tilde{b}_\Phi\|_{\tilde{L}^2_{t, \dot{\eta}(t)}(\mathcal{B}^{1, 0})}^2\right)\\
&\qquad+\left(4\alpha\|\hbar^{\frac12}e^\Psi\Delta_k^h\partial_y\tilde{u}_\Phi\|^2_{L^2(0, t; L_+^2)}+\frac{4C^2}\alpha\|\hbar^{\frac12}e^\Psi\Delta_k^h\partial_y\tilde{b}_\Phi\|^2_{L^2(0, t; L_+^2)}\right),
\end{split}
\end{equation}
and
\begin{equation}\label{449}
\begin{split}
\|\hbar^{\frac12}&e^\Psi\Delta_k^h\tilde{b}_\Phi\|_{L^\infty(0, t; L_+^2)}^2
+2c\lambda 2^k\int_0^t \dot{\eta}(t')\|\hbar^{\frac12}e^\Psi\Delta_k^h\tilde{b}_\Phi(t')\|_{L_+^2}^2\;dt'
+\|\hbar^{\frac12}e^\Psi\Delta_k^h\partial_y\tilde{b}_\Phi\|_{L^2(0, t; L_+^2)}^2\\
&\le\|\hbar^{\frac12}e^\Psi\Delta_k^h\tilde{b}_\Phi(0)\|_{L_+^2}^2
+\int_0^t \hbar '(t')\|e^\Psi\Delta_k^h\tilde{b}_\Phi\|_{L_+^2}^2\;dt'\\
&\qquad+Cd_k^2 2^{-k}\left(\|\hbar^{\frac12}e^\Psi\tilde{u}_\Phi\|_{\tilde{L}^2_{t, \dot{\eta}(t)}(\mathcal{B}^{1, 0})}^2+\|\hbar^{\frac12}e^\Psi\tilde{b}_\Phi\|_{\tilde{L}^2_{t, \dot{\eta}(t)}(\mathcal{B}^{1, 0})}^2\right).
\end{split}
\end{equation}
By taking square root of (\ref{448}) and (\ref{449}) and then multiplying the resulting inequalities by $2^{\frac k2}$ and finally summing over $k\in\mathbb{Z}$ respectively, we find that
\begin{equation}\label{450}
\begin{split}
\|\hbar^{\frac12}e^\Psi\tilde{u}_\Phi&\|_{\tilde{L}^\infty(0, t; \mathcal{B}^{\frac12, 0})}
+\sqrt{2c\lambda}\|\hbar^{\frac12}e^\Psi\tilde{u}_\Phi\|_{\tilde{L}^2_{t.\cdot{\eta}(t)}(\mathcal{B}^{1, 0})}
+\|\hbar^{\frac12}e^\Psi\partial_y\tilde{u}_\Phi\|_{\tilde{L}^2(0, t; \mathcal{B}^{\frac12,0})}\\
&\le\|\hbar^{\frac12}e^\Psi\tilde{u}_\Phi(0)\|_{\mathcal{B}^{\frac12, 0}}
+ \|\sqrt{\hbar '}e^\Psi\tilde{u}_\Phi\|_{\tilde{L}^2(0,t;\mathcal{B}^{\frac12,0})}\\
&\qquad+\sqrt{C}\left(\|\hbar^{\frac12}e^\Psi\tilde{u}_\Phi\|_{\tilde{L}^2_{t, \dot{\eta}(t)}(\mathcal{B}^{1, 0})}+\|\hbar^{\frac12}e^\Psi\tilde{b}_\Phi\|_{\tilde{L}^2_{t, \dot{\eta}(t)}(\mathcal{B}^{1, 0})}\right)\\
&\qquad+\left(2\sqrt{\alpha}\|\hbar^{\frac12}e^\Psi\partial_y\tilde{u}_\Phi\|_{\tilde{L}^2(0, t; \mathcal{B}^{\frac12,0})}+\frac {2C}{\sqrt{\alpha}}\|\hbar^{\frac12}e^\Psi\partial_y\tilde{b}_\Phi\|_{\tilde{L}^2(0, t; \mathcal{B}^{\frac12,0})}\right),
\end{split}
\end{equation}
and
\begin{equation}\label{451}
\begin{split}
\|\hbar^{\frac12}e^\Psi\tilde{b}_\Phi&\|_{\tilde{L}^\infty(0, t; \mathcal{B}^{\frac12, 0})}
+\sqrt{2c\lambda}\|\hbar^{\frac12}e^\Psi\tilde{b}_\Phi\|_{\tilde{L}^2_{t,\cdot{\eta}(t)}(\mathcal{B}^{1, 0})}
+\|\hbar^{\frac12}e^\Psi\partial_y\tilde{b}_\Phi\|_{\tilde{L}^2(0, t; \mathcal{B}^{\frac12,0})}\\
&\le\|\hbar^{\frac12}e^\Psi\tilde{b}_\Phi(0)\|_{\mathcal{B}^{\frac12, 0}}
+ \|\sqrt{\hbar '}e^\Psi\tilde{b}_\Phi\|_{\tilde{L}^2(0,t;\mathcal{B}^{\frac12,0})}\\
&\qquad+\sqrt{C}\left(\|\hbar^{\frac12}e^\Psi\tilde{u}_\Phi\|_{\tilde{L}^2_{t, \dot{\eta}(t)}(\mathcal{B}^{1, 0})}+\|\hbar^{\frac12}e^\Psi\tilde{b}_\Phi\|_{\tilde{L}^2_{t, \dot{\eta}(t)}(\mathcal{B}^{1, 0})}\right).
\end{split}
\end{equation}
Consider (\ref{450})+$K\times$(\ref{451}) with $K>1$ to be determined later and taking $\lambda$ to be large enough so that $2c\lambda\ge C$. In this way, we arrive at (\ref{49}). Thus,  the proof of Proposition \ref{P3.1} is done.

\section{Global Well-Posedness of Solutions}
In this section, we give another definition of the quantity $\theta(t)$ in (\ref{210}) to describe the evolution of the analytic band to the solution of (\ref{17}) after the moment of lower bound of lifespan of solutions, which is achieved in section 3.
\begin{equation}
\label{X1}
\begin{cases}
\dot\theta(t)=\langle t\rangle^{\frac14}\left(\|e^\Psi\partial_yG^s\|_{L_v^2}+f(t)\|e^\Psi\varrho '\|_{L_v^2}+\|e^\Psi\partial_y(G, H)_\Phi(t)\|_{\mathcal{B}^{\frac12, 0}}\right)\\
\theta|_{t=t_0}=0,
\end{cases}
\end{equation}
with $t_0=\frac1\varepsilon-1$. It should be pointed out that the lifespan obtained in section 3 is $\varepsilon^{-3/2}$, here and after $t_0=\frac1\varepsilon-1$ is enough for the analysis below.
The phase function $\Phi$ and the weighted function $\Psi$ are defined as in section 2. The goal of this section is to control the analytic radius and get the global well-posedness of solutions from $t_0$. Consequently, it suffices to estimate three terms on the right hand side of (\ref{X1}).

\subsection{Estimates of $G^s$}
The goal of this subsection is to give an estimate of $G^s$ defined by (\ref{27}).
Indeed, we have the following lemma.
\begin{lem}
\label{L4}
Let $G^s$ be defined as (\ref{27}). Then for any $t>t_0$ with $t_0=\frac 1\varepsilon-1$, one has
\[
\int_{t_0}^t \langle t'\rangle^{\frac14}\|e^\Psi\partial_yG^s(t')\|_{L_v^2}\;dt'\le CC_f\varepsilon^{\frac18}
\]
for some constant $C$, and $C_f$ is defined in (\ref{112}).
\end{lem}
\textbf{Proof of Lemma \ref{L4}}
By taking $L_v^2$ inner product of the equation of $G^s$ (\ref{27}) with $e^{2\Psi}G^s$, we have
\[
(\partial_tG^s |e^{2\Psi}G^s)_{L_v^2}-(\partial_y^2G^s |e^{2\Psi}G^s)_{L_v^2}+\langle t\rangle^{-1}\|e^\Psi G^s(t)\|_{L_v^2}^2=(\tilde{m} |e^{2\Psi}G^s)_{L_v^2}.
\]
First,
\[
(\partial_tG^s |e^{2\Psi}G^s)_{L_v^2}=\frac12\frac d{dt}\|e^\Psi G^s(t)\|_{L_v^2}^2-\int_{\mathbb{R}^+} e^{2\Psi}\partial_t\Psi |G^s|^2\;dy.
\]
Since $G^s|_{y=0}$, Using integration by parts and Young's inequality give that
\begin{align*}
-(\partial_y^2G^s |e^{2\Psi}G^s)_{L_v^2}&=\|e^\Psi\partial_yG^s(t)\|_{L_v^2}^2+2\int_{\mathbb{R}^+} e^{2\Psi}\partial_y\Psi\partial_yG^sG^s\;dy\\
&\ge\frac12 \|e^\Psi\partial_yG^s(t)\|_{L_v^2}^2-2\int_{\mathbb{R}^+} e^{2\Psi}(\partial_y\Psi)^2 |G^s|^2\;dy.
\end{align*}
Thanks to (\ref{214}), we arrive at
\begin{equation}\label{51}
\frac12\frac d{dt}\|e^\Psi G^s(t)\|_{L_v^2}^2+\frac12 \|e^\Psi\partial_yG^s(t)\|_{L_v^2}^2+\langle t\rangle^{-1}\|e^\Psi G^s(t)\|_{L_v^2}^2\le\|e^\Psi G^s(t)\|_{L_v^2}\|e^\Psi \tilde{m}(t)\|_{L_v^2}.
\end{equation}
In particular, we have
\[
\frac12\frac d{dt}\|e^\Psi G^s(t)\|_{L_v^2}^2+\langle t\rangle^{-1}\|e^\Psi G^s(t)\|_{L_v^2}^2\le\|e^\Psi G^s(t)\|_{L_v^2}\|e^\Psi \tilde{m}(t)\|_{L_v^2}.
\]
Dividing by $\|e^\Psi G^s(t)\|_{L_v^2}$ on both sides, we obtain
\[
\frac d{dt}\|e^\Psi G^s(t)\|_{L_v^2}+\langle t\rangle^{-1}\|e^\Psi G^s(t)\|_{L_v^2}\le\|e^\Psi \tilde{m}(t)\|_{L_v^2}.
\]
Moreover,  by multiplying the above inequality with $\langle t\rangle$ and integrating it over $[0, t]$, we get
\begin{equation}\label{52}
\|\langle t\rangle e^\Psi G^s(t, \cdot)\|_{L_v^2}\le\|\langle t'\rangle e^\Psi\tilde{m}\|_{L^1(0, t; L_v^2)},
\end{equation}
where we used the fact that $G^s|_{t=0}=0$. On the other hand, we deduce from (\ref{51}) and Young's inequality that
\begin{equation}\label{53}
\begin{split}
\frac d{dt}\|e^\Psi G^s(t)\|_{L_v^2}^2+\|e^\Psi\partial_yG^s(t)\|_{L_v^2}^2&+2\langle t\rangle^{-1}\|e^\Psi G^s(t)\|_{L_v^2}^2\\
&\le 2\langle t\rangle^{\frac12}\|e^\Psi\tilde{m}(t)\|_{L_v^2}\langle t\rangle^{-\frac12}\|e^\Psi G^s(t)\|_{L_v^2}\\
&\le \langle t\rangle \|e^\Psi\tilde{m}(t)\|_{L_v^2}^2+\langle t\rangle^{-1}\|e^\Psi G^s(t)\|_{L_v^2}^2.
\end{split}
\end{equation}
Multiplying the above inequality by $\langle t\rangle^2$ and then integrating it over $[t/2, t]$, we obtain
\begin{align*}
\int_{\frac t2}^t \|\langle t'\rangle e^\Psi\partial_y G^s(t')\|_{L_v^2}^2\;dt'&\le\|\langle t/2\rangle e^\Psi G^s(t/2)\|_{L_v^2}^2+\int_{\frac t2}^t \langle t'\rangle \|e^\Psi G^s(t')\|_{L_v^2}^2\;dt'+\int_{\frac t2}^t\langle t'\rangle^3 \|e^\Psi\tilde{m}(t)\|_{L_v^2}^2\;dt'\\
&\le(1+\ln 2)\|\langle t'\rangle e^\Psi G^s(t')\|_{L^\infty(0, t; L_v^2)}^2+
\|\langle t'\rangle^{\frac32}e^\Psi\tilde{m}\|_{L^2(0, t; L_v^2)}^2,
\end{align*}
which together with (\ref{52}) implies that
\begin{equation}\label{54}
\int_{\frac t2}^t \|\langle t'\rangle e^\Psi\partial_y G^s(t')\|_{L_v^2}^2\;dt'\le C\left(
\|\langle t'\rangle e^\Psi\tilde{m}\|_{L^1(0, t; L_v^2)}^2
+\|\langle t'\rangle^{\frac32}e^\Psi\tilde{m}\|_{L^2(0, t; L_v^2)}^2\right).
\end{equation}
Now, we are in a position to prove
\begin{equation}\label{55}
\int_{t_0}^t \langle t'\rangle^{\frac14}\|e^\Psi\partial_yG^s(t')\|_{L_v^2}\;dt'\le CC_f\varepsilon^{\frac18}.
\end{equation}
In view of (\ref{112}) and notice that $\tilde{m}$ is supported in $[0, 2]$, we have
\begin{align*}
\|\langle t\rangle\tilde{m}\|_{L^1(0, t; L_v^2)}&\lesssim\|yM\|_{L^1(0, t; L_v^2)}+\|\langle t\rangle m\|_{L^1(0, t; L_v^2)}\\
&\le C\int_0^\infty \langle t\rangle (|f(t)|+|f'(t)|)\;dt\\
&\le CC_f
\end{align*}
and
\begin{align*}
\|\langle t\rangle^{\frac32}\tilde{m}\|_{L^2(0, t; L_v^2)}&\lesssim\|\langle t\rangle^{\frac12}yM\|_{L^2(0, t; L_v^2)}+\|\langle t\rangle^{\frac32}m\|_{L^2(0, t; L_v^2)}\\
&\le C\left(\int_0^\infty \langle t\rangle^3 (f^2(t)+(f'(t))^2)\;dt\right)^{\frac12}\\
&\le CC_f
\end{align*}
for $C_f$ defined in (\ref{112}). From which and (\ref{54}), we get that
\[
\int_{\frac t2}^t \|\langle t'\rangle e^\Psi\partial_y G^s(t')\|_{L_v^2}^2\;dt'\le CC_f^2.
\]
Moreover, by integrating (\ref{53}) over $[0, t]$, we obtain
\[
\|e^\Psi \partial_yG^s(t)\|_{L^2(0, t; L_v^2)}^2\le CC_f^2.
\]
Then, for any $t>1$, there exists a unique integer $N_t$ so that $2^{N_t-1}<t\le 2^{N_t}$. Thus we have $\frac t2\le 2^{N_t-1}$ and
\begin{align*}
\int_{2^{N_t-1}}^t \langle t'\rangle^{\frac38} \|e^\Psi\partial_yG^s(t')\|_{L_v^2}\;dt'&\le
\left(\int_{2^{N_t-1}}^t \langle t'\rangle^{-\frac54}dt'\right)^{\frac12}
\left(\int_{t/2}^t (\langle t'\rangle \|e^\Psi\partial_yG^s(t')\|_{L_v^2})^2\;dt'\right)^{\frac12}\\
&\le C2^{-\frac{N_t}8}C_f
\end{align*}
In a similar way, for any $j\in [0, N_t-2]$, we have
\begin{align*}
\int_{2^j}^{2^{j+1}} \langle t'\rangle^{\frac38} \|e^\Psi\partial_yG^s(t')\|_{L_v^2}\;dt'&\le
\left(\int_{2^j}^{2^{j+1}} \langle t'\rangle^{-\frac54}dt'\right)^{\frac12}
\left(\int_{2^j}^{2^{j+1}} (\langle t'\rangle \|e^\Psi\partial_yG^s(t')\|_{L_v^2}^2\;dt'\right)^{\frac12}\\
&\le C2^{-\frac j8}C_f.
\end{align*}
Next, we derive from the above three inequalities that
\begin{align*}
&\int_{0}^t \langle t'\rangle^{\frac38} \|e^\Psi\partial_yG^s(t')\|_{L_v^2}\;dt'\le
2^{\frac38}\int_0^1\|e^\Psi\partial_yG^s(t')\|_{L_v^2}\;dt'\\
&\qquad+\sum\limits_{j=0}^{N_t-2}\int_{2^j}^{2^{j+1}} \langle t'\rangle^{\frac38} \|e^\Psi\partial_yG^s(t')\|_{L_v^2}\;dt'
+\int_{2^{N_t-1}}^t \langle t'\rangle^{\frac38} \|e^\Psi\partial_yG^s(t')\|_{L_v^2}\;dt'\\
&\le C(1+\sum\limits_{j=0}^{\infty}2^{-\frac j8})C_f\\
&\le CC_f.
\end{align*}
Finally, it comes out
\[
\langle t_0\rangle^{\frac18}\int_{t_0}^t \langle t'\rangle^{\frac14}\|e^\Psi\partial_yG^s(t')\|_{L_v^2}\;dt'\le\int_{t_0}^t \langle t'\rangle^{\frac38}\|e^\Psi\partial_yG^s(t')\|_{L_v^2}\;dt'\le CC_f.
\]
Then, we complete the proof of Lemma \ref{L4} by using $t_0=\frac{1}{\varepsilon}-1$.

We also have the following corollary.
\begin{cor}
Let $u^s$ be determined by (\ref{16}). Then for any $\gamma\in (0, 1)$, we have
\begin{equation}\label{57}
\|e^{\gamma\Psi}\partial_yu^s(t)\|_{L_v^2}\le \|e^{\Psi}\partial_yG^s(t)\|_{L_v^2}.
\end{equation}
\end{cor}
The proof is similar to Lemma \ref{L3.2} and we omit the detail here.

\subsection{Estimates of $(G, H)$}
In this subsection, we will present the estimates of $(G, H)$ defined in (\ref{29}).
\begin{lem}
\label{L5}
Let $G$ and $H$ be defined by (\ref{29}). Then for any $t>t_0$ with $t_0=\frac1\varepsilon-1$, it holds
\[
\int_{t_0}^t \langle t'\rangle^{\frac14}\|e^\Psi\partial_y(G, H)_\Phi(t')\|_{\mathcal{B}^{\frac12, 0}}\;dt'\le C\varepsilon^{\frac14}.
\]
\end{lem}
The proof of Lemma \ref{L5} relies on the following proposition.
\begin{prop}\label{P5}
Let $G$ and $H$ be determined by (\ref{29}). Then for any $t$, we have
\begin{align*}
\|\langle t'\rangle e^\Psi(G, H)_\Phi\|_{\tilde{L}^\infty(0, t;\mathcal{B}^{\frac12, 0})}+\|\langle t'\rangle e^\Psi\partial_y(G, H)_\Phi\|_{\tilde{L}^2(0, t; \mathcal{B}^{\frac12, 0})}
\le C\|e^\Psi e^{\delta |D_x|}(G_0, H_0)\|_{\mathcal{B}^{\frac12, 0}}.
\end{align*}
\end{prop}
We also admit Proposition \ref{P5} for the time being and prove Lemma \ref{L5}.

Indeed, for any $t>t_0$ and $t_0=\frac 1\varepsilon-1$, by H\"older's inequality, we have
\begin{align*}
\int_{t_0}^t\langle t'\rangle^{\frac14}\|e^\Psi\partial_y(G, H)_\Phi(t')\|_{\mathcal{B}^{\frac12, 0}}\;dt'&\le
\left(\int_{t_0}^t\langle t'\rangle^{-\frac32}\;dt'\right)^{\frac12}\left(\int_{t_0}^t\|\langle t'\rangle e^\Psi\partial_y(G, H)_\Phi(t')\|_{\mathcal{B}^{\frac12, 0}}^2\;dt'\right)^{\frac12}\\
&\le C\langle t_0\rangle^{-\frac14}\|e^\Psi e^{\delta |D_x|}(G_0, H_0)\|_{\mathcal{B}^{\frac12, 0}}\\
&\le C\varepsilon^{\frac14}.
\end{align*}
Thus, the proof of Lemma \ref{L5} is complete.  It remains to prove Proposition \ref{P5}.

\textbf{Proof of Proposition \ref{P5}.} First, in view of (\ref{211}), it comes from (\ref{29}) that
\begin{equation}\label{57}
\begin{split}
\hbar(t)&(e^\Psi\Delta_k^h(\partial_tG_\Phi-\partial_y^2G_\Phi)|e^\Psi\Delta_k^hG_\Phi)_{L_+^2}
+\lambda\dot{\theta}(t)\hbar(t)(e^\Psi |D_h|\Delta_k^hG_\Phi |e^\Psi\Delta_k^hG_\Phi)_{L_+^2}\\
&+{\langle t\rangle}^{-1}\hbar(t)(e^\Psi\Delta_k^hG_\Phi |e^\Psi\Delta_k^hG_\Phi)_{L_+^2}
+\hbar(t)(e^\Psi \Delta_k^h[v\partial_y(u^s+f(t)\varrho(y))]_\Phi |e^\Psi\Delta_k^hG_\Phi)_{L_+^2}\\
&+\hbar(t)(e^\Psi\Delta_k^h[(u+u^s+f(t)\varrho(y))\partial_xG]_\Phi |e^\Psi\Delta_k^hG_\Phi)_{L_+^2}
+\hbar(t)(e^\Psi\Delta_k^h[v\partial_yG]_\Phi |e^\Psi\Delta_k^hG_\Phi)_{L_+^2}\\
&-\hbar(t)(e^\Psi\Delta_k^h[(1+b)\partial_xH]_\Phi |e^\Psi\Delta_k^hG_\Phi)_{L_+^2}
-\hbar(t)(e^\Psi\Delta_k^h[g\partial_yH]_\Phi |e^\Psi\Delta_k^hG_\Phi)_{L_+^2}\\
&-\hbar(t)(\frac12{\langle t\rangle}^{-1}e^\Psi\Delta_k^h[v\partial_y(y\varphi)]_\Phi |e^\Psi\Delta_k^hG_\Phi)_{L_+^2}
+\hbar(t)(\frac12{\langle t\rangle}^{-1}e^\Psi\Delta_k^h[g\partial_y(y\psi)]_\Phi |e^\Psi\Delta_k^hG_\Phi)_{L_+^2}\\
&-\langle t\rangle^{-1}\hbar(t)(e^\Psi y\int_y^\infty \Delta_k^h[v\partial_y(u+u^s+f(t)\varrho(y))]_\Phi\;dy' |e^\Psi\Delta_k^hG_\Phi)_{L_+^2}\\
&+\langle t\rangle^{-1}\hbar(t)(e^\Psi y\int_y^\infty \Delta_k^h[g\partial_yb]_\Phi\;dy' |e^\Psi\Delta_k^hG_\Phi)_{L_+^2}=0
\end{split}
\end{equation}
and
\begin{equation}\label{58}
\begin{split}
\hbar(t)&(e^\Psi\Delta_k^h(\partial_tH_\Phi-\partial_y^2H_\Phi)  |e^\Psi\Delta_k^hH_\Phi)_{L_+^2}
+\lambda\dot{\theta}(t)\hbar(t)(e^\Psi |D_h|\Delta_k^hH_\Phi |e^\Psi\Delta_k^hH_\Phi)_{L_+^2}\\
&+{\langle t\rangle}^{-1}\hbar(t)(e^\Psi\Delta_k^hH_\Phi |e^\Psi\Delta_k^hH_\Phi)_{L_+^2}
-\hbar(t)(e^\Psi \Delta_k^h[g\partial_y(u^s+f(t)\varrho(y))]_\Phi |e^\Psi\Delta_k^hH_\Phi)_{L_+^2}\\
&+\hbar(t)(e^\Psi\Delta_k^h[(u+u^s+f(t)\varrho(y))\partial_xH]_\Phi |e^\Psi\Delta_k^hH_\Phi)_{L_+^2}
+\hbar(t)(e^\Psi\Delta_k^h[v\partial_yH]_\Phi |e^\Psi\Delta_k^hH_\Phi)_{L_+^2}\\
&-\hbar(t)(e^\Psi\Delta_k^h[(1+b)\partial_xG]_\Phi |e^\Psi\Delta_k^hH_\Phi)_{L_+^2}
-\hbar(t)(e^\Psi\Delta_k^h[g\partial_yG]_\Phi |e^\Psi\Delta_k^hH_\Phi)_{L_+^2}\\
&+\hbar(t)(\frac12{\langle t\rangle}^{-1}e^\Psi\Delta_k^h[g\partial_y(y\varphi)]_\Phi |e^\Psi\Delta_k^hH_\Phi)_{L_+^2}
-\hbar(t)(\frac12{\langle t\rangle}^{-1}e^\Psi\Delta_k^h[v\partial_y(y\psi)]_\Phi |e^\Psi\Delta_k^hH_\Phi)_{L_+^2}=0
\end{split}
\end{equation}
Let us integrate (\ref{57}) and (\ref{58}) over $[0,t]$ and then handle term by term in two resulting equalities.\\
$\bullet$ Estimate of $\int_{0}^t \hbar(t')(e^\Psi\Delta_k^h(\partial_tG_\Phi-\partial_y^2G_\Phi) |e^\Psi\Delta_k^hG_\Phi)_{L_+^2}\;dt^\prime$ and

\qquad\qquad$\int_{0}^t \hbar(t')(e^\Psi\Delta_k^h(\partial_tH_\Phi-\partial_y^2H_\Phi) |e^\Psi\Delta_k^hH_\Phi)_{L_+^2}\;dt^\prime$\\

Integration by parts leads to
\begin{equation}\label{59}
\begin{split}
&\int_{0}^t \hbar(t')(e^\Psi\Delta_k^h\partial_tG_\Phi |e^\Psi\Delta_k^hG_\Phi)_{L_+^2}\;dt^\prime\\
=&\frac12\|\hbar^{\frac12}e^\Psi\Delta_k^hG_\Phi(t)\|_{L_+^2}^2-\frac12\|\hbar^{\frac12}e^\Psi\Delta_k^hG_\Phi(0)\|_{L_+^2}^2\\
&-\frac12\int_{0}^t \hbar^\prime(t^\prime)\|e^\Psi\Delta_k^hG_\Phi(t^\prime)\|_{L_+^2}^2\;dt^\prime-\int_{0}^t\int_{\mathbb{R}_+^2} \hbar\partial_t\Psi |e^\Psi\Delta_k^hG_\Phi|_{L_+^2}^2\;dxdydt^\prime,
\end{split}
\end{equation}
and
\begin{equation}\label{510}
\begin{split}
&-\int_{0}^t \hbar(t')(e^\Psi\Delta_k^h\partial_y^2G_\Phi |e^\Psi\Delta_k^hG_\Phi)_{L_+^2}\;dt^\prime\\
=&\|\hbar^{\frac12}e^\Psi\Delta_k^h\partial_yG_\Phi(t)\|_{L^2(0, t;L_+^2)}^2+2\int_{0}^t\int_{\mathbb{R}_+^2} \hbar\partial_y\Psi e^{2\Psi}\Delta_k^hG_\Phi\Delta_k^h\partial_yG_\Phi\;dxdydt^\prime\\
\ge&\frac12\|\hbar^{\frac12}e^\Psi\Delta_k^h\partial_yG_\Phi(t)\|_{L^2(0, t;L_+^2)}^2
-2\int_{0}^t\int_{\mathbb{R}_+^2} \hbar(\partial_y\Psi)^2 |e^{\Psi}\Delta_k^hG_\Phi |^2\;dxdydt^\prime.
\end{split}
\end{equation}
Due to (\ref{214}), we conclude that
\begin{equation}\label{511}
\begin{split}
\int_{0}^t &\hbar(t')(e^\Psi\Delta_k^h(\partial_tG_\Phi-\partial_y^2G_\Phi) |e^\Psi\Delta_k^hG_\Phi)_{L_+^2}\;dt^\prime\\
&\ge\frac12\left(\|\hbar^{\frac12}e^\Psi\Delta_k^hG_\Phi(t)\|_{L_+^2}^2
-\|\hbar^{\frac12}e^\Psi\Delta_k^hG_\Phi(0)\|_{L_+^2}^2\right.\\
&\left.\quad-\int_{0}^t \hbar^\prime(t^\prime)\|e^\Psi\Delta_k^hG_\Phi(t^\prime)\|_{L_+^2}^2\;dt^\prime
+\|\hbar^{\frac12}e^\Psi\Delta_k^h\partial_yG_\Phi(t)\|_{L^2(0, t;L_+^2)}^2\right).
\end{split}
\end{equation}
With a similar argument, we also get
\begin{equation}\label{512}
\begin{split}
\int_{0}^t &\hbar(t')(e^\Psi\Delta_k^h(\partial_tH_\Phi-\partial_y^2H_\Phi) |e^\Psi\Delta_k^hH_\Phi)_{L_+^2}\;dt^\prime\\
&\ge\frac12\left(\|\hbar^{\frac12}e^\Psi\Delta_k^hH_\Phi(t)\|_{L_+^2}^2
-\|\hbar^{\frac12}e^\Psi\Delta_k^hH_\Phi(0)\|_{L_+^2}^2\right.\\
&\left.\quad-\int_{0}^t \hbar^\prime(t^\prime)\|e^\Psi\Delta_k^hH_\Phi(t^\prime)\|_{L_+^2}^2\;dt^\prime
+\|\hbar^{\frac12}e^\Psi\Delta_k^h\partial_yH_\Phi(t)\|_{L^2(0, t;L_+^2)}^2\right).
\end{split}
\end{equation}
$\bullet$ Estimate of $\int_0^t\lambda\dot{\theta}(t')\hbar(t')(e^\Psi |D_h|\Delta_k^h G_\Phi |e^\Psi\Delta_k^h G_\Phi)_{L_+^2}\;dt^\prime$ and

\qquad\qquad$\int_0^t\lambda\dot{\theta}(t')\hbar(t')(e^\Psi |D_h|\Delta_k^h H_\Phi |e^\Psi\Delta_k^h H_\Phi)_{L_+^2}\;dt^\prime$\\

From Lemma \ref{L2}, we get that
\begin{equation}\label{513}
\int_{0}^t\lambda\dot{\theta}(t')\hbar(t')(e^\Psi |D_h|\Delta_k^hG_\Phi |e^\Psi\Delta_k^hG_\Phi)_{L_+^2}\;dt^\prime
\ge c\lambda 2^k\int_{0}^t\dot{\theta}(t')\|\hbar^{\frac12}e^\Psi\Delta_k^hG_\Phi\|_{L_+^2}^2\;dt^\prime,
\end{equation}
and
\begin{equation}\label{514}
\int_{0}^t\lambda\dot{\theta}(t')\hbar(t')(e^\Psi |D_h|\Delta_k^hH_\Phi |e^\Psi\Delta_k^hH_\Phi)_{L_+^2}\;dt^\prime
\ge c\lambda 2^k\int_{0}^t\dot{\theta}(t')\|\hbar^{\frac12}e^\Psi\Delta_k^hH_\Phi\|_{L_+^2}^2\;dt^\prime.
\end{equation}
$\bullet$ Estimate of $\int_{0}^t \hbar(t')(e^\Psi\Delta_k^h[(u+u^s+f(t')\varrho(y))\partial_xG]_\Phi |e^\Psi\Delta_k^hG_\Phi)_{L_+^2}\;dt^\prime$ and

\qquad\qquad$\int_{0}^t \hbar(t')(e^\Psi\Delta_k^h[(u+u^s+f(t')\varrho(y))\partial_xH]_\Phi |e^\Psi\Delta_k^hH_\Phi)_{L_+^2}\;dt^\prime$

Integration by parts yields
\begin{equation}\label{515}
\int_{0}^t \hbar(t')(e^\Psi\Delta_k^h[(u^s+f(t')\varrho(y))\partial_xG]_\Phi |e^\Psi\Delta_k^hG_\Phi)_{L_+^2}\;dt^\prime=0
\end{equation}
Then, we apply Bony's decomposition on $u\partial_xG$ to obtain
\[
u\partial_xG=T_u^h\partial_xG+T_{\partial_xG}^hu+R^h(u, \partial_xG).
\]
By virtue of (\ref{213}) and  the support properties of the Fourier transform of the terms in $T_u^h\partial_xG$, we have
\begin{align*}
&\int_{0}^t \hbar(t')|(e^\Psi\Delta_k^h[T_u^h\partial_xG]_\Phi |e^\Psi\Delta_k^hG_\Phi)_{L_+^2}|\;dt^\prime\\
\lesssim&\sum\limits_{|k'-k|\le4}\int_{0}^t \|S_{k'-1}^hu
_\Phi(t')\|_{L_+^\infty}\|\hbar^{\frac12}e^\Psi\Delta_{k'}^h\partial_xG_\Phi(t')\|_{L_+^2}\|\hbar^{\frac12}e^\Psi\Delta_k^hG_\Phi(t')\|_{L_+^2}\;dt'.
\end{align*}
Notice that
\begin{equation}\label{516}
\begin{split}
\|\Delta_k^hu_\Phi(t')\|_{L_v^\infty(L_h^2)}&\lesssim\|\Delta_k^h(\int_0^y\partial_yu_\Phi(t')\;dy')\|_{L_v^\infty(L_h^2)}\\
&\lesssim\|\Delta_k^h\partial_yu_\Phi(t')\|_{L_v^1(L_h^2)}\\
&\lesssim\|e^{-\frac{\Psi}4}\|_{L_v^2} \|e^{\frac{\Psi}4}\Delta_k^h\partial_yu_\Phi(t')\|_{L_+^2}\\
&\lesssim d_k(t')2^{-\frac k2}\langle t'\rangle^{\frac14}\|e^\Psi\partial_yG_\Phi(t')\|_{\mathcal{B}^{\frac12, 0}},
\end{split}
\end{equation}
in the last inequality, Lemma \ref{L3.2} is used. Then, by Lemma 2.1,
\begin{align*}
\|S_{k'-1}^hu_\Phi(t')\|_{L_+^\infty}&\lesssim\sum\limits_{k\le k'-2} 2^{\frac k2}\|\Delta_k^hu_\Phi(t')|_{L_v^\infty(L_h^2)}\\
&\lesssim\langle t'\rangle^{\frac14}\|e^\Psi\partial_yG_\Phi(t')\|_{\mathcal{B}^{\frac12, 0}}\lesssim\dot{\theta}(t').
\end{align*}
Thus, we get from H\"older's inequality that
\begin{equation}\label{517}
\begin{split}
&\int_0^t \hbar(t')|(e^\Psi\Delta_k^h[T_u^h\partial_xG]_\Phi |e^\Psi\Delta_k^hG_\Phi)_{L_+^2}|\;dt^\prime\\
\lesssim& 2^{k'}\left(\int_{0}^t\dot{\theta}(t')\|\hbar^{\frac12}e^\Psi\Delta_{k'}^hG_\Phi(t')\|_{L_+^2}^2\;dt'\right)^{\frac12}
\left(\int_{0}^t\dot{\theta}(t')\|\hbar^{\frac12}e^\Psi\Delta_{k}^hG_\Phi(t')\|_{L_+^2}^2\;dt'\right)^{\frac12}\\
\lesssim& d_k^22^{-k}\|\hbar^{\frac12}e^\Psi G_\Phi\|^2_{\tilde{L}_{t, \dot{\theta}(t)}^2(\mathcal{B}^{1,0})}.
\end{split}
\end{equation}
Similarly, by H\"older's inequality again, we have
\begin{align*}
&\int_{0}^t \hbar(t)|(e^\Psi\Delta_k^h[T_{\partial_xG}^hu]_\Phi |e^\Psi\Delta_k^hG_\Phi)_{L_+^2}|\;dt^\prime\\
\lesssim&\sum\limits_{|k'-k|\le4}\int_{0}^t \|\hbar^{\frac12}e^\Psi S_{k'-1}^h\partial_xG
_\Phi(t')\|_{L_v^2(L_h^\infty)}
\|\Delta_{k'}^hu_\Phi(t')\|_{L_v^\infty(L_h^2)}
\|\hbar^{\frac12}e^\Psi\Delta_k^hG_\Phi(t')\|_{L_+^2}\;dt'\\
\lesssim&\sum\limits_{|k'-k|\le4} d_{k'}2^{-\frac {k'}2}\left(\int_{0}^t\dot{\theta}(t')\|\hbar^{\frac12}e^\Psi S_{k'-1}^h\partial_xG_\Phi(t')\|_{L_v^2(L_h^\infty)}^2\;dt'\right)^{\frac12}
\left(\int_{0}^t\dot{\theta}(t')\|\hbar^{\frac12}e^\Psi\Delta_{k}^hG_\Phi(t')\|_{L_+^2}^2\;dt'\right)^{\frac12}\\
\end{align*}
where (\ref{516}) and the definition of $\dot{\theta}(t)$ are used in the last inequality. By Lemma \ref{L2},
\begin{align*}
\left(\int_0^t\dot{\theta}(t')\|\hbar^{\frac12}e^\Psi S_{k'-1}^h\partial_xG_\Phi(t')\|_{L_v^2(L_h^\infty)}^2\;dt'\right)^{\frac12}
&\lesssim\sum\limits_{j\le k'-2}2^{\frac{3j}2}\left(\int_0^t\dot{\theta}(t')\|\hbar^{\frac12}e^\Psi\Delta_{j}^hG_\Phi(t')\|_{L_+^2}^2\;dt'\right)^{\frac12}\\
&\lesssim 2^{\frac{k'}2}\|\hbar^{\frac12}e^\Psi G_\Phi\|_{\tilde{L}_{t, \dot{\theta}(t)}^2(\mathcal{B}^{1,0})}.
\end{align*}
Hence, it comes out
\begin{equation}\label{518}
\int_{0}^t \hbar(t')|(e^\Psi\Delta_k^h[T_{\partial_xG}^hu]_\Phi |e^\Psi\Delta_k^hG_\Phi)_{L_+^2}|\;dt^\prime
\lesssim d_k^22^{-k}\|\hbar^{\frac12}e^{\Psi} G_\Phi\|^2_{\tilde{L}_{t, \dot{\theta}(t)}^2(\mathcal{B}^{1,0})}.
\end{equation}
For the last term, by Lemma \ref{L2} and (\ref{516}) again, we have
\begin{align*}
&\int_{0}^t \hbar(t')|(e^\Psi\Delta_k^h[R^h(u, \partial_xG)]_\Phi |e^\Psi\Delta_k^hG_\Phi)_{L_+^2}|\;dt^\prime\\
\lesssim&\sum\limits_{k'\ge k-3}\int_{0}^t \|\tilde{\Delta}_{k'}^hu_\Phi(t')\|_{L_+^\infty}
\|\hbar^{\frac12}e^\Psi\Delta_{k'}^h\partial_xG_\Phi(t')\|_{L_+^2}
\|\hbar^{\frac12}e^\Psi\Delta_k^hG_\Phi(t')\|_{L_+^2}\;dt'\\
\lesssim&\sum\limits_{k'\ge k-3}2^{k'}\int_0^t \dot{\theta}(t')
\|\hbar^{\frac12}e^\Psi\Delta_{k'}^hG_\Phi(t')\|_{L_+^2}
\|\hbar^{\frac12}e^\Psi\Delta_k^hG_\Phi(t')\|_{L_+^2}\;dt',
\end{align*}
from which we obtain that
\begin{equation}\label{519}
\int_{0}^t \hbar(t')|(e^\Psi\Delta_k^h[R^h(u, \partial_xG)]_\Phi |e^\Psi\Delta_k^hG_\Phi)_{L_+^2}|\;dt^\prime\lesssim d_k^22^{-k}\|\hbar^{\frac12}e^\Psi G_\Phi\|^2_{\tilde{L}_{t, \dot{\theta}(t)}^2(\mathcal{B}^{1,0})}.
\end{equation}
Combining (\ref{515}) and (\ref{517})-(\ref{519}) together gives
\begin{equation}\label{520}
\int_{0}^t \hbar(t')|(e^\Psi\Delta_k^h[(u+u^s+f(t')\varrho(y))\partial_xG]_\Phi |e^\Psi\Delta_k^hG_\Phi)_{L_+^2}|\;dt^\prime\lesssim d_k^22^{-k}\|\hbar^{\frac12}e^\Psi G_\Phi\|^2_{\tilde{L}_{t, \dot{\theta}(t)}^2(\mathcal{B}^{1,0})}.
\end{equation}
Along the same line, we also have
\begin{equation}\label{521}
\int_{0}^t \hbar(t')|(e^\Psi\Delta_k^h[(u+u^s+f(t')\varrho(y))\partial_xH]_\Phi |e^\Psi\Delta_k^hH_\Phi)_{L_+^2}|\;dt^\prime\lesssim d_k^22^{-k}\|\hbar^{\frac12}e^\Psi H_\Phi\|^2_{\tilde{L}_{t, \dot{\theta}(t)}^2(\mathcal{B}^{1,0})}.
\end{equation}
$\bullet$Estimate of $\int_0^t\hbar(t')(e^\Psi \Delta_k^h[(1+b)\partial_xH]_\Phi |e^\Psi\Delta_k^hG_\Phi)_{L_+^2}\;dt^\prime$ and

 \qquad\qquad$\int_0^t\hbar(t')(e^\Psi \Delta_k^h[(1+b)\partial_xG]_\Phi |e^\Psi\Delta_k^hH_\Phi)_{L_+^2}\;dt^\prime$\\

Indeed,  the estimates of these two terms are similar as the above estimates. We only present these estimates and omit the detail here.
\begin{equation}\label{522}
\begin{split}
&\int_{0}^t \hbar(t')|(e^\Psi\Delta_k^h[(1+b)\partial_xH]_\Phi |e^\Psi\Delta_k^hG_\Phi)_{L_+^2}|\;dt^\prime\\
\lesssim& d_k^22^{-k}\left(\|\hbar^{\frac12}e^\Psi G_\Phi\|_{\tilde{L}_{t, \dot{\theta}(t)}^2(\mathcal{B}^{1,0})}^2+\|\hbar^{\frac12}e^\Psi H_\Phi\|_{\tilde{L}_{t, \dot{\theta}(t)}^2(\mathcal{B}^{1,0})}^2\right),
\end{split}
\end{equation}
and
\begin{equation}\label{523}
\begin{split}
&\int_{0}^t \hbar(t')|(e^\Psi\Delta_k^h[(1+b)\partial_xG]_\Phi |e^\Psi\Delta_k^hH_\Phi)_{L_+^2}|\;dt^\prime\\
\lesssim& d_k^22^{-k}\left(\|\hbar^{\frac12}e^\Psi G_\Phi\|_{\tilde{L}_{t, \dot{\theta}(t)}^2(\mathcal{B}^{1,0})}^2+\|\hbar^{\frac12}e^\Psi H_\Phi\|_{\tilde{L}_{t, \dot{\theta}(t)}^2(\mathcal{B}^{1,0})}^2\right).
\end{split}
\end{equation}
$\bullet$ Estimate of $\int_{0}^t\hbar(t')(e^\Psi \Delta_k^h[v\partial_yG]_\Phi |e^\Psi\Delta_k^hG_\Phi)_{L_+^2}\;dt'$, $\int_{0}^t\hbar(t')(e^\Psi \Delta_k^h[v\partial_yH]_\Phi |e^\Psi\Delta_k^hH_\Phi)_{L_+^2}\;dt'$\\

 \qquad and$\int_{0}^t\hbar(t')(e^\Psi \Delta_k^h[g\partial_yH]_\Phi |e^\Psi\Delta_k^hG_\Phi)_{L_+^2}\;dt'$, $\int_{0}^t\hbar(t')(e^\Psi \Delta_k^h[g\partial_yG]_\Phi |e^\Psi\Delta_k^hH_\Phi)_{L_+^2}\;dt'$

Use Bony's decomposition for $v\partial_yG$ as before and get
\[
v\partial_yG=T_v^h\partial_yG+T_{\partial_yG}^hv+R^h(v, \partial_yG).
\]
In view of (\ref{213}) and the support properties of the Fourier transform of the terms in $T_v^h\partial_yG$, we write
\begin{align*}
&\int_{0}^t \hbar(t')|(e^\Psi\Delta_k^h[T_v^h\partial_yG]_\Phi |e^\Psi\Delta_k^hG_\Phi)_{L_+^2}|\;dt^\prime\\
\lesssim&\sum\limits_{|k'-k|\le4}\int_{0}^t \|\hbar^{\frac12}S_{k'-1}^hv_\Phi(t')\|_{L_+^\infty}
\|e^\Psi\Delta_{k'}^h\partial_yG_\Phi(t')\|_{L_+^2}\|\hbar^{\frac12}e^\Psi\Delta_k^hG_\Phi(t')\|_{L_+^2}\;dt'\\
\lesssim&\sum\limits_{|k'-k|\le4}d_{k'} 2^{-\frac{k'}2}\int_0^t \|\hbar^{\frac12}S_{k'-1}^h(\int_0^y\partial_xu_\Phi(t')\;dy')\|_{L_+^\infty}
\|e^\Psi\partial_yG_\Phi(t')\|_{\mathcal{B}^{\frac12, 0}}\|\hbar^{\frac12}e^\Psi\Delta_k^hG_\Phi(t')\|_{L_+^2}\;dt'.
\end{align*}
Besides, by a similar argument in (\ref{516}), we have
\begin{align*}
\|S_{k'-1}^h(\int_0^y\partial_xu_\Phi(t')\;dy')\|_{L_+^\infty}&\lesssim\sum\limits_{k\le k'-2}2^{\frac32 k}\|\Delta_{k}^hu_\Phi (t')\|_{L_v^1(L_h^2)}\\
&\lesssim\sum\limits_{k\le k'-2}2^{\frac32 k}\|e^{-\frac{\Psi}4}\|_{L_v^2}\|e^{\frac{\Psi}4}\Delta_{k}^hu_\Phi (t')\|_{L_+^2}\\
&\lesssim\sum\limits_{k\le k'-2}2^{\frac32 k}\langle t'\rangle^{\frac14}\|e^\Psi\Delta_{k}^hG_\Phi (t')\|_{L_+^2},
\end{align*}
where in the last inequality, Lemma \ref{L3.2} is used. Then, we obtain
\begin{equation}\label{524}
\int_{0}^t \hbar(t')|(e^\Psi\Delta_k^h[T_v^h\partial_xG]_\Phi |e^\Psi\Delta_k^hG_\Phi)_{L_+^2}|\;dt^\prime
\lesssim d_k^22^{-k}\|\hbar^{\frac12}e^\Psi G_\Phi\|_{\tilde{L}_{t, \dot{\theta}(t)}^2(\mathcal{B}^{1,0})}^2.
\end{equation}
By the same way, we also have
\begin{align*}
&\int_{0}^t \hbar(t')|(e^\Psi\Delta_k^h[T_{\partial_yG}^hv]_\Phi |e^\Psi\Delta_k^hG_\Phi)_{L_+^2}|\;dt^\prime\\
\lesssim&\sum\limits_{|k'-k|\le4}\int_{0}^t
\|e^\Psi S_{k'-1}^h\partial_yG_\Phi(t')\|_{L_v^2(L_h^\infty)}
\|\hbar^{\frac12}\Delta_{k'}^h(\int_0^y \partial_yv_\Phi(t')\;dy')\|_{L_v^\infty(L_h^2)}
\|\hbar^{\frac12}e^\Psi\Delta_k^hG_\Phi(t')\|_{L_+^2}\;dt'.
\end{align*}
Whereas in view of Lemmas \ref{L2} and \ref{L3.2}, we get
\begin{equation*}
\|e^\Psi S_{k'-1}^h\partial_yG_\Phi(t')\|_{L_v^2(L_h^\infty)}\lesssim\sum\limits_{k\le k'-2}2^{\frac{k}2}\|e^\Psi\Delta_{k}^h\partial_yG_\Phi (t')\|_{L_+^2}\lesssim\|e^{\Psi} \partial_yG_\Phi (t')\|_{\mathcal{B}^{\frac12, 0}},
\end{equation*}
and
\begin{equation}\label{525}
\begin{split}
\|\hbar^{\frac12}\Delta_{k'}^h(\int_0^y \partial_yv_\Phi(t')\;dy')\|_{L_v^\infty(L_h^2)}&\lesssim
2^{k'}\|\hbar^{\frac12}\Delta_{k'}^h(\int_0^y u_\Phi(t')\;dy')\|_{L_v^\infty(L_h^2)}\\
&\lesssim
2^{k'}\|\hbar^{\frac12}\Delta_{k'}^h u_\Phi(t')\|_{L_v^1(L_h^2)}\\
&\lesssim
2^{k'}\|e^{-\frac{\Psi}4}(t')\|_{L_v^2}\|\hbar^{\frac12}e^{\frac\Psi 4}\Delta_{k'}^h u_\Phi(t')\|_{L_+^2}\\
&\lesssim 2^{k'}\langle t'\rangle ^{\frac14}\|\hbar^{\frac12}e^\Psi\Delta_{k'}^h G_\Phi(t')\|_{L_+^2}.
\end{split}
\end{equation}
Hence, we arrive at the following estimate by collecting the above inequalities.
\begin{equation}\label{526}
\int_{0}^t \hbar(t')|(e^\Psi\Delta_k^h[T_{\partial_yG}^hv]_\Phi |e^\Psi\Delta_k^hG_\Phi)_{L_+^2}|\;dt^\prime
\lesssim d_k^22^{-k}\|\hbar^{\frac12}e^\Psi G_\Phi\|_{\tilde{L}_{t, \dot{\theta}(t)}^2(\mathcal{B}^{1,0})}^2.
\end{equation}
For the last term,
\begin{align*}
&\int_{0}^t \hbar(t')|(e^\Psi\Delta_k^h[R^h(v, \partial_yG)]_\Phi |e^\Psi\Delta_k^hG_\Phi)_{L_+^2}|\;dt^\prime\\
\lesssim&\sum\limits_{k'\ge k-3}\int_{0}^t
\|\hbar^{\frac12}\Delta_{k'}^h(\int_0^y\partial_yv_\Phi(t')\;dy')\|_{L_v^\infty(L_h^2)}
\|e^\Psi\tilde{\Delta}_{k'}^h \partial_yG_\Phi(t')\|_{L_v^2(L_h^\infty)}
\|\hbar^{\frac12}e^\Psi\Delta_k^hG_\Phi(t')\|_{L_+^2}\;dt'.
\end{align*}
Using (\ref{525}), we gain
\begin{equation}\label{527}
\int_{0}^t \hbar(t')|(e^\Psi\Delta_k^h[R^h(v, \partial_yG)]_\Phi |e^\Psi\Delta_k^hG_\Phi)_{L_+^2}|\;dt^\prime
\lesssim d_k^22^{-k}\|\hbar^{\frac12}e^\Psi G_\Phi\|_{\tilde{L}_{t, \dot{\theta}(t)}^2(\mathcal{B}^{1,0})}^2.
\end{equation}
As a consequence,
\begin{equation}\label{528}
\int_{0}^t\hbar(t')|(e^\Psi \Delta_k^h[v\partial_yG]_\Phi |e^\Psi\Delta_k^hG_\Phi)_{L_+^2}|\;dt^\prime
\lesssim d_k^22^{-k}\|\hbar^{\frac12}e^\Psi G_\Phi\|_{\tilde{L}_{t, \dot{\theta}(t)}^2(\mathcal{B}^{1,0})}^2.
\end{equation}
Along the same line, it holds that
\begin{equation}\label{529}
\int_{0}^t\hbar(t')|(e^\Psi \Delta_k^h[g\partial_yG]_\Phi |e^\Psi\Delta_k^hH_\Phi)_{L_+^2}|\;dt^\prime
\lesssim d_k^22^{-k}\|\hbar^{\frac12}e^\Psi H_\Phi\|_{\tilde{L}_{t, \dot{\theta}(t)}^2(\mathcal{B}^{1,0})}^2,
\end{equation}
\begin{equation}\label{530}
\begin{split}
&\int_{0}^t\hbar(t')|(e^\Psi\Delta_k^h[g\partial_yH]_\Phi |e^\Psi\Delta_k^hG_\Phi)_{L_+^2}|\;dt^\prime\\
\lesssim& d_k^22^{-k}\left(\|\hbar^{\frac12}e^\Psi G_\Phi\|_{\tilde{L}_{t, \dot{\theta}(t)}^2(\mathcal{B}^{1,0})}^2+\|\hbar^{\frac12}e^\Psi H_\Phi\|_{\tilde{L}_{t, \dot{\theta}(t)}^2(\mathcal{B}^{1,0})}^2\right),
\end{split}
\end{equation}
and
\begin{equation}\label{531}
\begin{split}
&\int_{0}^t\hbar(t')|(e^\Psi\Delta_k^h[v\partial_yH]_\Phi |e^\Psi\Delta_k^hH_\Phi)_{L_+^2}|\;dt^\prime\\
\lesssim& d_k^22^{-k}\left(\|\hbar^{\frac12}e^\Psi G_\Phi\|_{\tilde{L}_{t, \dot{\theta}(t)}^2(\mathcal{B}^{1,0})}^2+\|\hbar^{\frac12}e^\Psi H_\Phi\|_{\tilde{L}_{t, \dot{\theta}(t)}^2(\mathcal{B}^{1,0})}^2\right).
\end{split}
\end{equation}
$\bullet$ Estimate of $\int_{0}^t\hbar(t')(e^\Psi \Delta_k^h[v\partial_y(u^s+f(t')\varrho(y))]_\Phi |e^\Psi\Delta_k^hG_\Phi)_{L_+^2}\;dt'$ and

\qquad\qquad$\int_{0}^t\hbar(t')(e^\Psi \Delta_k^h[g\partial_y(u^s+f(t')\varrho(y))]_\Phi |e^\Psi\Delta_k^hH_\Phi)_{L_+^2}\;dt'$\\

First,
\begin{align*}
&\int_{0}^t\hbar(t')|(e^\Psi \Delta_k^h[v\partial_y(u^s+f(t')\varrho(y))]_\Phi |e^\Psi\Delta_k^hG_\Phi)_{L_+^2}|\;dt'\\
\lesssim& \int_{0}^t (\|e^{\frac34\Psi} \partial_yu^s\|_{L_v^2}+f(t')\|e^{\frac34\Psi}\varrho'\|_{L_v^2})
\|\hbar^{\frac12}e^{\frac\Psi 4}\Delta_k^hv_\Phi(t')\|_{L_v^\infty(L_h^2)}
\|\hbar^{\frac12}e^\Psi\Delta_k^hG_\Phi\|_{L_+^2}\;dt',
\end{align*}
combining which and (\ref{525}) imply that
\begin{equation}\label{532}
\int_{0}^t\hbar(t')|(e^\Psi \Delta_k^h[v\partial_y(u^s+f(t')\varrho(y))]_\Phi |e^\Psi\Delta_k^hG_\Phi)_{L_+^2}|\;dt'
\lesssim d_k^22^{-k}\|\hbar^{\frac12}e^\Psi G_\Phi\|_{\tilde{L}_{t, \dot{\theta}(t)}^2(\mathcal{B}^{1,0})}^2,
\end{equation}
where the following two inequalities are used.
\begin{align*}
\|e^{\frac34\Psi}\partial_yu^s\|_{L_v^2}\lesssim \|e^{\Psi}\partial_yG^s\|_{L_v^2}\lesssim \dot{\theta}(t)\langle t'\rangle^{-1/4},
\end{align*}
and
\begin{align*}
f(t')\|e^{\frac34\Psi}\varrho'\|_{L_v^2}\lesssim \|e^{\Psi}\varrho'\|_{L_v^2}\lesssim \dot{\theta}(t)\langle t'\rangle^{-1/4}.
\end{align*}
Similarly,
\begin{equation}\label{533}
\int_{0}^t\hbar(t')|(e^\Psi \Delta_k^h[g\partial_y(u^s+f(t')\varrho(y))]_\Phi |e^\Psi\Delta_k^hH_\Phi)_{L_+^2}|\;dt'
\lesssim d_k^22^{-k}\|\hbar^{\frac12}e^\Psi H_\Phi\|_{\tilde{L}_{t, \dot{\theta}(t)}^2(\mathcal{B}^{1,0})}^2.
\end{equation}
$\bullet$ Estimate of $\int_{0}^t\hbar(t')\langle t'\rangle^{-1}(e^\Psi \Delta_k^h[v\partial_y(y\varphi)]_\Phi |e^\Psi\Delta_k^hG_\Phi)_{L_+^2}\;dt'$,\\
$\int_{0}^t\hbar(t')\langle t'\rangle^{-1}(e^\Psi \Delta_k^h[g\partial_y(y\varphi)]_\Phi |e^\Psi\Delta_k^hH_\Phi)_{L_+^2}\;dt'$, $\int_{0}^t\hbar(t')\langle t'\rangle^{-1}(e^\Psi \Delta_k^h[g\partial_y(y\psi)]_\Phi |e^\Psi\Delta_k^hG_\Phi)_{L_+^2}\;dt'$,
\qquad and $\int_{0}^t\hbar(t')\langle t'\rangle^{-1}(e^\Psi \Delta_k^h[v\partial_y(y\psi)]_\Phi |e^\Psi\Delta_k^hH_\Phi)_{L_+^2}\;dt'$

 Applying Bony's decomposition on $v\partial_y(y\varphi)$ yields
 \[
 v\partial_y(y\varphi)=T_v^h\partial_y(y\varphi)+T_{\partial_y(y\varphi)}v+R^h(v, \partial_y(y\varphi)).
 \]
From (\ref{213}) and the support properties of the terms in $T_v^h\partial_y(y\varphi)$, we deduce that
\begin{align*}
&\int_{0}^t\hbar(t')\langle t'\rangle^{-1}|(e^\Psi \Delta_k^h[T_v^h\partial_y(y\varphi)]_\Phi |e^\Psi\Delta_k^hG_\Phi)_{L_+^2}|\;dt'\\
\lesssim&\sum\limits_{|k'-k|\le4}\int_{0}^t \langle t'\rangle^{-1}\|\hbar^{\frac12}e^{\frac\Psi4}S_{k'-1}^hv_\Phi(t')\|_{L_+^\infty}
\|e^{\frac{3\Psi}4}\Delta_{k'}^h\partial_y(y\varphi_\Phi)(t')\|_{L_+^2}
\|\hbar^{\frac12}e^\Psi\Delta_k^hG_\Phi(t')\|_{L_+^2}\;dt'\\
\lesssim&\sum\limits_{|k'-k|\le4}\int_{0}^t \|\hbar^{\frac12}e^{\frac\Psi4}S_{k'-1}^hv_\Phi(t')\|_{L_+^\infty}
\|e^\Psi\Delta_{k'}^h\partial_yG_\Phi(t')\|_{L_+^2}
\|\hbar^{\frac12}e^\Psi\Delta_k^hG_\Phi(t')\|_{L_+^2}\;dt',
\end{align*}
here Lemma \ref{L3.2} is used in the last inequality. Then, by a similar argument as in (\ref{516}), we achieve
\begin{align*}
\|e^{\frac\Psi4}S_{k'-1}^hv_\Phi(t')\|_{L_+^\infty}&\lesssim\sum\limits_{k\le k'-2}2^{\frac32 k}\|e^{\frac\Psi 4}\Delta_{k}^hu_\Phi (t')\|_{L_v^1(L_h^2)}\\
&\lesssim\sum\limits_{k\le k'-2}2^{\frac32 k}\|e^{-\frac{\Psi}4} (t')\|_{L_v^2}\|e^{\frac{\Psi}2}\Delta_{k}^hu_\Phi (t')\|_{L_+^2}\\
&\lesssim\sum\limits_{k\le k'-2}2^{\frac32 k}\langle t'\rangle^{\frac14}\|e^\Psi\Delta_{k}^hG_\Phi (t')\|_{L_+^2}.
\end{align*}
As a result, we deduce that
\begin{equation}\label{534}
\int_{0}^t\hbar(t')\langle t'\rangle^{-1}|(e^\Psi \Delta_k^h[T_v^h\partial_y(y\varphi)]_\Phi |e^\Psi\Delta_k^hG_\Phi)_{L_+^2}|\;dt'\lesssim d_k^22^{-k}\|\hbar^{\frac12}e^\Psi G_\Phi\|_{\tilde{L}_{t, \dot{\theta}(t)}^2(\mathcal{B}^{1,0})}^2.
\end{equation}
Similarly,
\begin{align*}
&\int_{0}^t\hbar(t')\langle t'\rangle^{-1}|(e^\Psi \Delta_k^h[T_{\partial_y(y\varphi)}^hv]_\Phi |e^\Psi\Delta_k^hG_\Phi)_{L_+^2}|\;dt'\\
\lesssim&\sum\limits_{|k'-k|\le4}\int_{0}^t \langle t'\rangle^{-1}
\|e^{\frac{3\Psi}4}S_{k'-1}^h\partial_y(y\varphi_\Phi)(t')\|_{L_v^2(L_h^\infty)}
\|\hbar^{\frac12}e^{\frac\Psi4}\Delta_{k'}^hv_\Phi(t')\|_{L_v^\infty(L_h^2)}
\|\hbar^{\frac12}e^\Psi\Delta_k^hG_\Phi(t')\|_{L_+^2}\;dt'.
\end{align*}
On the other hand, we have
\begin{align*}
\langle t'\rangle^{-1}\|e^{\frac{3\Psi}4}S_{k'-1}^h\partial_y(y\varphi_\Phi)(t')\|_{L_v^2(L_h^\infty)}
&\lesssim\sum\limits_{j\le k'-2}2^{\frac{j}2}\langle t'\rangle^{-1}\|e^{\frac{3\Psi}4}\Delta_{j}^h\partial_y(y\varphi_\Phi)(t')\|_{L_+^2}\\
&\lesssim\sum\limits_{j\le k'-2} 2^{\frac {j}2}\|e^{\Psi}\Delta_{j}^h\partial_yG_\Phi(t')\|_{L_+^2}\\
&\lesssim \|e^{\Psi}\partial_yG_\Phi(t')\|_{\mathcal{B}^{\frac12, 0}}\lesssim \langle t'\rangle^{-\frac14}\dot{\theta}(t'),
\end{align*}
and
\begin{align*}
\|e^{\frac\Psi4}\Delta_{k'}^hv_\Phi(t')\|_{L_v^\infty(L_h^2)}&\lesssim 2^{k'}\|e^{\frac\Psi4}\Delta_{k'}^hu_\Phi (t')\|_{L_v^1(L_h^2)}\\
&\lesssim 2^{ k'}\|e^{-\frac{\Psi}4} (t')\|_{L_v^2}\|e^{\frac{\Psi}2}\Delta_{k'}^hu_\Phi (t')\|_{L_+^2}\\
&\lesssim 2^{ k'}\langle t'\rangle^{\frac14}\|e^\Psi\Delta_{k'}^hG_\Phi (t')\|_{L_+^2},
\end{align*}
Then, we reach
\begin{equation}\label{535}
\int_{0}^t\hbar(t')\langle t'\rangle^{-1}|(e^\Psi \Delta_k^h[T_{\partial_y(y\varphi)}^hv]_\Phi |e^\Psi\Delta_k^hG_\Phi)_{L_+^2}|\;dt'
\lesssim d_k^22^{-k}\|\hbar^{\frac12}e^\Psi G_\Phi\|_{\tilde{L}_{t, \dot{\theta}(t)}^2(\mathcal{B}^{1,0})}^2.
\end{equation}
Finally, again thanks to the support properties of the Fourier transform of the terms in $R^h(v, \partial_y(y\varphi))$, we get
\begin{align*}
&\int_{0}^t\hbar(t')\langle t'\rangle^{-1}|(e^\Psi \Delta_k^h[R^h(v, \partial_y(y\varphi))]_\Phi |e^\Psi\Delta_k^hG_\Phi)_{L_+^2}|\;dt'\\
\lesssim&\sum\limits_{k'\ge k-3}\int_{0}^t \langle t'\rangle^{-1}\|\hbar^{\frac12}e^{\frac\Psi4}\Delta_{k'}^hv_\Phi(t')\|_{L_+^\infty}
\|e^{\frac{3\Psi}4}\tilde{\Delta}_{k'}^h\partial_y(y\varphi_\Phi)(t')\|_{L_+^2}
\|\hbar^{\frac12}e^\Psi\Delta_k^hG_\Phi(t')\|_{L_+^2}\;dt'\\
\lesssim&\sum\limits_{k'\ge k-3}\int_{0}^t \|\hbar^{\frac12}e^{\frac\Psi4}\Delta_{k'}^hv_\Phi(t')\|_{L_+^\infty}
\|e^\Psi\Delta_{k'}^h\partial_yG_\Phi(t')\|_{L_+^2}
\|\hbar^{\frac12}e^\Psi\Delta_k^hG_\Phi(t')\|_{L_+^2}\;dt'.
\end{align*}
Then, it follows from a similar argument in (\ref{525}) that
\begin{equation}\label{536}
\int_{0}^t\hbar(t')\langle t'\rangle^{-1}|(e^\Psi \Delta_k^h[R^h(v, \partial_y(y\varphi))]_\Phi |e^\Psi\Delta_k^hG_\Phi)_{L_+^2}|\;dt'
\lesssim d_k^22^{-k}\|\hbar^{\frac12}e^\Psi G_\Phi\|_{\tilde{L}_{t, \dot{\theta}(t)}^2(\mathcal{B}^{1,0})}^2.
\end{equation}
By summing up the above estimates, we conclude that
\begin{equation}\label{537}
\int_{0}^t\hbar(t')\langle t'\rangle^{-1}|(e^\Psi \Delta_k^h[v\partial_y(y\varphi)]_\Phi |e^\Psi\Delta_k^hG_\Phi)_{L_+^2}|\;dt'
\lesssim d_k^22^{-k}\|\hbar^{\frac12}e^\Psi G_\Phi\|_{\tilde{L}_{t, \dot{\theta}(t)}^2(\mathcal{B}^{1,0})}^2.
\end{equation}
By the same arguments, the following estimates also hold.
\begin{equation}\label{538}
\int_{0}^t\hbar(t')\langle t'\rangle^{-1}|(e^\Psi \Delta_k^h[g\partial_y(y\varphi)]_\Phi |e^\Psi\Delta_k^hH_\Phi)_{L_+^2}|\;dt'
\lesssim d_k^22^{-k}\|\hbar^{\frac12}e^\Psi H_\Phi\|_{\tilde{L}_{t, \dot{\theta}(t)}^2(\mathcal{B}^{1,0})}^2,
\end{equation}
\begin{equation}\label{539}
\begin{split}
&\int_{0}^t\hbar(t')\langle t'\rangle^{-1}|(e^\Psi \Delta_k^h[g\partial_y(y\psi)]_\Phi |e^\Psi\Delta_k^hG_\Phi)_{L_+^2}|\;dt'\\
\lesssim& d_k^22^{-k}\left(\|\hbar^{\frac12}e^\Psi G_\Phi\|_{\tilde{L}_{t, \dot{\theta}(t)}^2(\mathcal{B}^{1,0})}^2+\|\hbar^{\frac12}e^\Psi H_\Phi\|_{\tilde{L}_{t, \dot{\theta}(t)}^2(\mathcal{B}^{1,0})}^2\right),
\end{split}
\end{equation}
and
\begin{equation}\label{540}
\begin{split}
&\int_{0}^t\hbar(t')\langle t'\rangle^{-1}|(e^\Psi \Delta_k^h[v\partial_y(y\psi)]_\Phi |e^\Psi\Delta_k^hH_\Phi)_{L_+^2}|\;dt'\\
\lesssim& d_k^22^{-k}\left(\|\hbar^{\frac12}e^\Psi G_\Phi\|_{\tilde{L}_{t, \dot{\theta}(t)}^2(\mathcal{B}^{1,0})}^2+\|\hbar^{\frac12}e^\Psi H_\Phi\|_{\tilde{L}_{t, \dot{\theta}(t)}^2(\mathcal{B}^{1,0})}^2\right).
\end{split}
\end{equation}
$\bullet$ Estimate of $\int_{0}^t\hbar(t')\langle t'\rangle^{-1}(e^\Psi y\int_y^\infty \Delta_k^h[v\partial_y(u+u^s+f(t')\varrho (y))]_\Phi\;dy' |e^\Psi\Delta_k^hG_\Phi)_{L_+^2}\;dt'$ and

\qquad\qquad $\int_{0}^t\hbar(t')\langle t'\rangle^{-1}(e^\Psi y\int_y^\infty \Delta_k^h[g\partial_yb]_\Phi\;dy' |e^\Psi\Delta_k^hG_\Phi)_{L_+^2}\;dt'$

Since
\begin{align*}
&\int_{0}^t\hbar(t')\langle t'\rangle^{-1}|(e^\Psi y\int_y^\infty \Delta_k^h[\partial_y(u^s+f(t')\varrho (y))\cdot v]_\Phi\;dy' |e^\Psi\Delta_k^hG_\Phi)_{L_+^2}|\;dt'\\
\lesssim&\int_{0}^t\langle t'\rangle^{-1}\|e^{-\frac\Psi 4}y\|_{L_v^2}
(\|e^{{\frac34}\Psi} \partial_yu^s\|_{L_v^2}+f(t')\|e^{\frac34\Psi}\varrho'\|_{L_v^2})
\|\hbar^{\frac12}e^{\frac\Psi2}\Delta_k^hv_\Phi(t')\|_{L_+^2}
\|\hbar^{\frac12}e^\Psi\Delta_k^hG_\Phi \|_{L_+^2}\;dt'\\
\lesssim&\int_{0}^t \langle t'\rangle^{-\frac14}(\|e^{{\frac34}\Psi} \partial_yu^s\|_{L_v^2}+f(t')\|e^\Psi\varrho'\|_{L_v^2})
\|\hbar^{\frac12}e^{\frac\Psi2}\Delta_k^hv_\Phi(t')\|_{L_+^2}
\|\hbar^{\frac12}e^\Psi\Delta_k^hG_\Phi \|_{L_+^2}\;dt'
\end{align*}
and
\begin{equation}\label{541}
\begin{split}
\|e^{\frac\Psi2}\Delta_k^hv_\Phi(t)\|_{L_+^2}&\lesssim 2^k\|e^{\frac{\Psi}2}\int_y^\infty \Delta_k^hu_\Phi\;dy'\|_{L_+^2}\\
&\lesssim 2^k\|e^{-\frac{\Psi}8}\int_y^\infty e^{-\frac\Psi8}\times e^{\frac{3\Psi}4}\Delta_k^hu_\Phi\;dy'\|_{L_+^2}\\
&\lesssim 2^k\|e^{-\frac{\Psi}8}\|_{L_v^2}^2\|e^{\frac{3\Psi}4} \Delta_k^hu_\Phi(t)\|_{L_+^2}\\
&\lesssim 2^k \langle t\rangle^{\frac12}\|e^{\Psi} \Delta_k^hG_\Phi(t)\|_{L_+^2},
\end{split}
\end{equation}
we deduce that
\begin{equation}\label{542}
\begin{split}
&\int_{0}^t\hbar(t')\langle t'\rangle^{-1}|(e^\Psi y\int_y^\infty \Delta_k^h[\partial_y(u^s+f(t')\varrho (y))\cdot v]_\Phi\;dy' |e^\Psi\Delta_k^hG_\Phi)_{L_+^2}|\;dt'\\
\lesssim& d_k^22^{-k}\|\hbar^{\frac12}e^\Psi G_\Phi\|_{\tilde{L}_{t, \dot{\theta}(t)}^2(\mathcal{B}^{1,0})}^2.
\end{split}
\end{equation}
Again, applying Bony's decomposition on $v\partial_yu$ yields
\[
v\partial_yu=T_v^h\partial_yu+T_{\partial_yu}^hv+R^h(v, \partial_yu).
\]
Due to (\ref{213}) and the support properties of the Fourier transform of the terms in $T_v^h\partial_yu$, we get
\begin{align*}
&\int_{0}^t\hbar(t')\langle t'\rangle^{-1}|(e^\Psi y\int_y^\infty \Delta_k^h[T_v^h\partial_yu]_\Phi\;dy' |e^\Psi\Delta_k^hG_\Phi)_{L_+^2}|\;dt'\\
\lesssim&\sum\limits_{|k'-k|\le4}\int_{0}^t\langle t'\rangle^{-1}\|e^{-\frac\Psi 4}y\|_{L_v^2}
\|\hbar^{\frac12}e^{\frac\Psi2}S_{k'-1}^hv_\Phi(t')\|_{L_v^2(L_h^\infty)}
\|e^{\frac{3\Psi}4}\Delta_{k'}^h\partial_yu_\Phi(t')\|_{L_+^2}
\|\hbar^{\frac12}e^\Psi\Delta_k^hG_\Phi \|_{L_+^2}\;dt'\\
\lesssim&\sum\limits_{|k'-k|\le4}\int_{0}^t \langle t'\rangle^{-\frac14}
\|\hbar^{\frac12}e^{\frac\Psi2}S_{k'-1}^hv_\Phi(t')\|_{L_v^2(L_h^\infty)}
\|e^{\Psi}\Delta_{k'}^h\partial_yG_\Phi\|_{L_+^2}
\|\hbar^{\frac12}e^\Psi\Delta_k^hG_\Phi \|_{L_+^2}\;dt'.
\end{align*}
In view of (\ref{541}),  we have
\[
 \|\hbar^{\frac12}e^{\frac\Psi2}S_{k'-1}^hv_\Phi(t)\|_{L_v^2(L_h^\infty)}\lesssim\sum\limits_{k\le k'-2} 2^{\frac{3k}2}\langle t\rangle^{\frac12}\|\hbar^{\frac12}e^{\Psi} \Delta_k^hG_\Phi(t)\|_{L_+^2}.
\]
As a result, it follows that
\begin{equation}\label{543}
\int_{0}^t\hbar(t')\langle t'\rangle^{-1}|(e^\Psi y\int_y^\infty \Delta_k^h[T_v^h\partial_yu]_\Phi\;dy' |e^\Psi\Delta_k^hG_\Phi)_{L_+^2}|\;dt'\lesssim d_k^22^{-k}\|\hbar^{\frac12}e^\Psi G_\Phi\|_{\tilde{L}_{t, \dot{\theta}(t)}^2(\mathcal{B}^{1,0})}^2.
\end{equation}
By a similar procedure, we find
\begin{align*}
&\int_{0}^t\hbar(t')\langle t'\rangle^{-1}|(e^\Psi y\int_y^\infty \Delta_k^h[T_{\partial_yu}^hv]_\Phi\;dy' |e^\Psi\Delta_k^hG_\Phi)_{L_+^2}|\;dt'\\
\lesssim&\sum\limits_{|k'-k|\le4}\int_{0}^t\langle t'\rangle^{-1}\|e^{-\frac\Psi 4}y\|_{L_v^2}
\|e^{\frac{3\Psi}4}S_{k'-1}^h\partial_yu_\Phi(t')\|_{L_+^2}
\|\hbar^{\frac12}e^{\frac\Psi2}\Delta_{k'}^hv_\Phi(t')\|_{L_v^2(L_h^\infty)}
\|\hbar^{\frac12}e^\Psi\Delta_k^hG_\Phi \|_{L_+^2}\;dt'\\
\lesssim&\sum\limits_{|k'-k|\le4}\int_{0}^t \langle t'\rangle^{-\frac14}
\|e^{\frac{3\Psi}4}S_{k'-1}^h\partial_yu_\Phi(t')\|_{L_+^2}
\|\hbar^{\frac12}e^{\frac\Psi2}\Delta_{k'}^hv_\Phi(t')\|_{L_v^2(L_h^\infty)}
\|\hbar^{\frac12}e^\Psi\Delta_k^hG_\Phi \|_{L_+^2}\;dt'.
\end{align*}
Next, a similar argument yields
\begin{equation}\label{544}
\int_{0}^t\hbar(t')\langle t'\rangle^{-1}|(e^\Psi y\int_y^\infty \Delta_k^h[T_{\partial_yu}^hv]_\Phi\;dy' |e^\Psi\Delta_k^hG_\Phi)_{L_+^2}|\;dt'\lesssim d_k^22^{-k}\|\hbar^{\frac12}e^\Psi G_\Phi\|_{\tilde{L}_{t, \dot{\theta}(t)}^2(\mathcal{B}^{1,0})}^2.
\end{equation}
Finally, thanks to (\ref{213}) and the support properties of the Fourier transform of the terms in $R^h(v, \partial_yu)$, we have
\begin{align*}
&\int_{0}^t\hbar(t')\langle t'\rangle^{-1}|(e^\Psi y\int_y^\infty \Delta_k^h[R^h(v, \partial_yu)]_\Phi\;dy' |e^\Psi\Delta_k^hG_\Phi)_{L_+^2}|\;dt'\\
\lesssim&\sum\limits_{k'\ge k-3}\int_{0}^t\langle t'\rangle^{-1}\|e^{-\frac\Psi 4}y\|_{L_v^2}
\|\hbar^{\frac12}e^{\frac\Psi2}\tilde{\Delta}_{k'}^hv_\Phi(t')\|_{L_v^2(L_h^\infty)}
\|e^{\frac{3\Psi}4}\Delta_{k'}^h\partial_yu_\Phi(t')\|_{L_+^2}
\|\hbar^{\frac12}e^\Psi\Delta_k^hG_\Phi \|_{L_+^2}\;dt'\\
\lesssim&\sum\limits_{k'\ge k-3}\int_{0}^t \langle t'\rangle^{-\frac14}
\|\hbar^{\frac12}e^{\frac\Psi2}\tilde{\Delta}_{k'}^hv_\Phi(t')\|_{L_v^2(L_h^\infty)}
\|e^{\Psi}\Delta_{k'}^h\partial_yG_\Phi\|_{L_+^2}
\|\hbar^{\frac12}e^\Psi\Delta_k^hG_\Phi \|_{L_+^2}\;dt',
\end{align*}
from which and a similar derivation in (\ref{541}), we obtain
\begin{equation}\label{545}
\int_{0}^t\hbar(t')\langle t'\rangle^{-1}|(e^\Psi y\int_y^\infty \Delta_k^h[R^h(v, \partial_yu)]_\Phi\;dy' |e^\Psi\Delta_k^hG_\Phi)_{L_+^2}|\;dt'\lesssim d_k^22^{-k}\|\hbar^{\frac12}e^\Psi G_\Phi\|_{\tilde{L}_{t, \dot{\theta}(t)}^2(\mathcal{B}^{1,0})}^2.
\end{equation}
Thus, we conclude that
\begin{equation}\label{546}
\begin{split}
&\int_{0}^t\hbar(t')\langle t'\rangle^{-1}|(e^\Psi y\int_y^\infty \Delta_k^h[v\partial_y(u+u^s+f(t')\varrho(y))]_\Phi\;dy' |e^\Psi\Delta_k^hG_\Phi)_{L_+^2}|\;dt'\\
\lesssim& d_k^22^{-k}\|\hbar^{\frac12}e^\Psi G_\Phi\|_{\tilde{L}_{t, \dot{\theta}(t)}^2(\mathcal{B}^{1,0})}^2.
\end{split}
\end{equation}
Similarly, we also have
\begin{equation}\label{547}
\begin{split}
&\int_{0}^t\hbar(t')\langle t'\rangle^{-1}|(e^\Psi y\int_y^\infty \Delta_k^h[g\partial_yb]_\Phi\;dy' |e^\Psi\Delta_k^hG_\Phi)_{L_+^2}|\;dt'\\
\lesssim& d_k^22^{-k}\left(\|\hbar^{\frac12}e^\Psi G_\Phi\|_{\tilde{L}_{t, \dot{\theta}(t)}^2(\mathcal{B}^{1,0})}^2+\|\hbar^{\frac12}e^\Psi H_\Phi\|_{\tilde{L}_{t, \dot{\theta}(t)}^2(\mathcal{B}^{1,0})}^2\right).
\end{split}
\end{equation}
Now,  integrating (\ref{57}) and (\ref{58}) over $[0, t]$ respectively, and inserting the above estimates into the resulting equalities, we achieve
\begin{equation}\label{548}
\begin{split}
&\|\hbar^{\frac12}e^\Psi\Delta_k^hG_\Phi\|_{L^\infty(0, t; L_+^2)}^2+2\|\langle t'\rangle^{-\frac12}\hbar^{\frac12}e^\Psi\Delta_k^hG_\Phi\|_{L^2(0, t; L_+^2)}^2
+\|\hbar^{\frac12}e^\Psi\Delta_k^h\partial_yG_\Phi\|_{L^2(0, t; L_+^2)}^2\\
&\qquad+2c\lambda 2^k\int_{0}^t \dot{\theta}(t')\|\hbar^{\frac12}e^\Psi\Delta_k^hG_\Phi(t')\|_{L_+^2}^2\;dt'\\
\le&\|\hbar^{\frac12}e^\Psi\Delta_k^hG_\Phi(0)\|_{L_+^2}^2
+\int_{0}^t \hbar '(t')\|e^\Psi\Delta_k^hG_\Phi\|_{L_+^2}^2\;dt'\\
&\qquad+Cd_k^2 2^{-k}\left(\|\hbar^{\frac12}e^\Psi G_\Phi\|_{\tilde{L}^2_{t, \dot{\theta}(t)}(\mathcal{B}^{1, 0})}^2+\|\hbar^{\frac12}e^\Psi H_\Phi\|_{\tilde{L}^2_{t, \dot{\theta}(t)}(\mathcal{B}^{1, 0})}^2\right)
\end{split}
\end{equation}
and
\begin{equation}\label{549}
\begin{split}
&\|\hbar^{\frac12}e^\Psi\Delta_k^hH_\Phi\|_{L^\infty(0, t; L_+^2)}^2+2\|\langle t'\rangle^{-\frac12}\hbar^{\frac12}e^\Psi\Delta_k^hH_\Phi\|_{L^2(0, t; L_+^2)}^2
+\|\hbar^{\frac12}e^\Psi\Delta_k^h\partial_yH_\Phi\|_{L^2(0, t; L_+^2)}^2\\
&\qquad+2c\lambda 2^k\int_{0}^t \dot{\theta}(t')\|\hbar^{\frac12}e^\Psi\Delta_k^hH_\Phi(t')\|_{L_+^2}^2\;dt'\\
\le&\|\hbar^{\frac12}e^\Psi\Delta_k^hH_\Phi(0)\|_{L_+^2}^2
+\int_{0}^t \hbar '(t')\|e^\Psi\Delta_k^hH_\Phi\|_{L_+^2}^2\;dt'\\
&\qquad+Cd_k^2 2^{-k}\left(\|\hbar^{\frac12}e^\Psi G_\Phi\|_{\tilde{L}^2_{t, \dot{\theta}(t)}(\mathcal{B}^{1, 0})}^2+\|\hbar^{\frac12}e^\Psi H_\Phi\|_{\tilde{L}^2_{t, \dot{\theta}(t)}(\mathcal{B}^{1, 0})}^2\right)
\end{split}
\end{equation}
By taking square root of (\ref{548}), (\ref{549}), then multiplying the resulting inequalities by $2^{\frac k2}$ and finally summing them over $k\in\mathbb{Z}$ respectively, we conclude that
\begin{equation}\label{550}
\begin{split}
&\|\hbar^{\frac12}e^\Psi G_\Phi\|_{\tilde{L}^\infty(0, t; \mathcal{B}^{\frac12, 0})}+\sqrt{2}\|\langle t'\rangle^{-\frac12}\hbar^{\frac12}e^\Psi G_\Phi\|_{\tilde{L}^2(0, t; \mathcal{B}^{\frac12,0})}
+\|\hbar^{\frac12}e^\Psi\partial_yG_\Phi\|_{\tilde{L}^2(0, t; \mathcal{B}^{\frac12, 0})}\\
&\qquad+\sqrt{2c\lambda} \|\hbar^{\frac12}e^\Psi G_\Phi\|_{\tilde{L}^2_{t; \dot{\theta}(t)}(\mathcal{B}^{1, 0})}\\
\le&\|\hbar^{\frac12}e^\Psi G_\Phi(0)\|_{\mathcal{B}^{\frac12,0}}
+\|\sqrt{\hbar '}e^\Psi G_\Phi\|_{\tilde{L}^2(0, t;\mathcal{B}^{\frac12, 0})}\\
&\qquad+\sqrt{C}\left(\|\hbar^{\frac12}e^\Psi G_\Phi\|_{\tilde{L}^2_{t, \dot{\theta}(t)}(\mathcal{B}^{1, 0})}+\|\hbar^{\frac12}e^\Psi H_\Phi\|_{\tilde{L}^2_{t, \dot{\theta}(t)}(\mathcal{B}^{1, 0})}\right),
\end{split}
\end{equation}
and
\begin{equation}\label{551}
\begin{split}
&\|\hbar^{\frac12}e^\Psi H_\Phi\|_{\tilde{L}^\infty(0, t; \mathcal{B}^{\frac12, 0})}+\sqrt{2}\|\langle t'\rangle^{-\frac12}\hbar^{\frac12}e^\Psi H_\Phi\|_{\tilde{L}^2(0, t; \mathcal{B}^{\frac12,0})}
+\|\hbar^{\frac12}e^\Psi\partial_yH_\Phi\|_{\tilde{L}^2(t; \mathcal{B}^{\frac12, 0})}\\
&\qquad+\sqrt{2c\lambda} \|\hbar^{\frac12}e^\Psi H_\Phi\|_{\tilde{L}^2_{t; \dot{\theta}(t)}(\mathcal{B}^{1, 0})}\\
\le&\|\hbar^{\frac12}e^\Psi H_\Phi(0)\|_{\mathcal{B}^{\frac12,0}}
+\|\sqrt{\hbar '}e^\Psi H_\Phi\|_{\tilde{L}^2(0, t;\mathcal{B}^{\frac12, 0})}\\
&\qquad+\sqrt{C}\left(\|\hbar^{\frac12}e^\Psi G_\Phi\|_{\tilde{L}^2_{t, \dot{\theta}(t)}(\mathcal{B}^{1, 0})}+\|\hbar^{\frac12}e^\Psi H_\Phi\|_{\tilde{L}^2_{t, \dot{\theta}(t)}(\mathcal{B}^{1, 0})}\right).
\end{split}
\end{equation}
Adding (\ref{550}) and (\ref{551}) together and taking $\lambda$ large enough, such that $2c\lambda\ge C$, we deduce
\begin{equation}\label{552}
\begin{split}
&\|\hbar^{\frac12}e^\Psi(G, H)_\Phi\|_{\tilde{L}^\infty(0, t; \mathcal{B}^{\frac12, 0})}
+\sqrt{2}\|\langle t'\rangle^{-\frac12}\hbar^{\frac12}e^\Psi(G, H)_\Phi\|_{\tilde{L}^2(0, t; \mathcal{B}^{\frac12,0})}\\
&+\|\hbar^{\frac12}e^\Psi\partial_y(G, H )_\Phi\|_{\tilde{L}^2(0, t; \mathcal{B}^{\frac12,0})}\\
\le& \|\hbar^{\frac12}e^\Psi(G, H)_\Phi(0)\|_{\mathcal{B}^{\frac12, 0}}
+ \|\sqrt{\hbar '}e^\Psi(G, H)_\Phi\|_{\tilde{L}^2(0, t;\mathcal{B}^{\frac12,0})}.
\end{split}
\end{equation}
Taking $\hbar(t)=\langle t\rangle^{2}$ in the above inequality, we find for any $t$, that
\begin{equation}\label{553}
\|\langle t'\rangle e^\Psi(G, H)_\Phi\|_{\tilde{L}^\infty(0, t; \mathcal{B}^{\frac12, 0})}+\|\langle t'\rangle e^\Psi\partial_y(G, H )_\Phi\|_{\tilde{L}^2(0, t; \mathcal{B}^{\frac12,0})} \le
\|e^{\Psi}e^{\delta |D_x|}(G_0, H_0)\|_{\mathcal{B}^{\frac12, 0}}.
\end{equation}
This completes the proof of Proposition \ref{P5}.

\subsection{Proof of main theorem \ref{T1.1}}
At the end of this section, we present the proof of main theorem.\\
{\bf \emph{Proof of Theorem \ref{T1.1}}:} The general strategy to prove the existence part in Theorem \ref{T1.1} is to construct an appropriate approximate solution sequence, and then to derive the uniform estimates for such an approximate solution sequence, finally to pass limit in the approximate problem by some compactness theorems. For simplicity, here we only present the uniform in time a priori estimates for the solutions to (\ref{17}) in the analytic framework, which we have achieved in the previous parts of this section. And the uniqueness of solutions can be proved by the same arguments as in \cite{XF}.

Let $(u, b)$ and $(\varphi, \psi)$ be smooth enough solutions to (\ref{17}) and (\ref{111}) respectively on $[0, t), t>t_0$. Let $(G, H )$ be defined by (\ref{29}). Then for any $t<T^*$ with $T^*$ being defined by (\ref{212}), we deduce that for $t_0=\frac1\varepsilon-1$,
\[
\theta(t)\le\int_{t_0}^t\langle t'\rangle^{\frac14}(\|e^\Psi\partial_yG^s\|_{L_v^2}+\|e^\Psi\partial_y(G, H)_\Phi(t')\|_{\mathcal{B}^{\frac12, 0}})\;dt'+\int_{t_0}^t\|e^\Psi\varrho'\|_{L_v^2}\langle t'\rangle^{\frac14}f(t')\;dt'.
\]
Notice that $\varrho'$ is supported in $[1, 2]$, and
\[
\int_{t_0}^t \langle t'\rangle^{\frac14}f(t')\;dt'\le\varepsilon^{\frac34}\int_{t_0}^t\langle t'\rangle |f(t')|\;dt'\le\varepsilon^{\frac34} C_f,
\]
from which, Lemmas \ref{L4} and \ref{L5}, we have
\begin{equation}\label{556}
\theta(t)\le C\varepsilon^{\frac18}
\end{equation}
for $t>t_0$, $\varepsilon<1$ and $C_f$ is determined in (\ref{112}).

In particular, if we take $\varepsilon$ small enough, such that
\[
C\varepsilon^{\frac18}\le\frac\delta{2\lambda}.
\]
Then we infer from (\ref{556}) that
\[
\sup\limits_{t\in [t_0, T^*)}\theta(t)\le\frac\delta{2\lambda}.
\]
In view of (\ref{212}), we get by a continuity argument that $T^*=\infty$. This completes the global existence part of Theorem 1.1. Moreover, (\ref{410}) holds for $t=\infty$, which together with (\ref{45}), (\ref{47}) and (\ref{48}) implies (\ref{114}).

\vspace{.1in}

\noindent{\bf Acknowledgments:}
The research of F. Xie is supported by National Natural Science Foundation of China No.11831003 and Shanghai Science and Technology Innovation Action Plan No. 20JC1413000.

\end{document}